\documentclass[11pt]{article}
\usepackage{verbatim}
\usepackage{amsfonts,amsmath,amsthm}
\def\LaTeX{L\kern -.36em\raise .3ex\hbox{\sc a}\kern -.15em T\kern -.1667em%
\lower .7ex\hbox{E}\kern -.125em X}

\usepackage{amssymb}
\usepackage{float}
\usepackage{epsfig}
\usepackage{amsmath}
\usepackage[english]{babel}
\usepackage{a4}
\usepackage{amstext}
\usepackage{mathrsfs}
\usepackage{amsthm}
\usepackage{fancyhdr}
\usepackage{latexsym}
\usepackage{amsmath}
\usepackage{amsfonts}
\usepackage{amssymb}
\usepackage{fancyhdr}

\usepackage[usenames]{color}
\definecolor{red}{rgb}{1.0,0.0,0.0}

\definecolor{blu}{rgb}{0.0,0.0,1.0}

\def\to{\longrightarrow}
\def \a {{\`{a} }}

\def\norm{{\| \kern -.05em | }}

\newtheorem{theorem}{Theorem}[section]
\newtheorem{proposition}[theorem]{Proposition}
\newtheorem{lemma}[theorem]{Lemma}

\theoremstyle{definition}
\newtheorem{definition}[theorem]{Definition}
\newtheorem{remark}[theorem]{Remark}

\newtheorem{hypothesis}[theorem]{Hypothesis}

\author{Salvatore Federico \and Ben Goldys \and Fausto Gozzi}

\begin{document}

\title{
HJB Equations for the Optimal Control of Differential Equations with
Delays and State Constraints, II: Optimal Feedbacks and
Approximations\footnote{This work was partially supported by an
Australian Research Council Discovery Project}}
\author{
Salvatore Federico\footnote{\thinspace Salvatore Federico,
Dipartimento di Scienze Economiche ed Aziendali, Facolt\`{a} di
Economia, Libera Universit\`{a} internazionale degli studi sociali
``Guido Carli'', viale Romania 32, 00197 Roma, Italy. Email:
\texttt{sfederico at luiss.it}.} \and Ben Goldys\footnote{\thinspace
Ben Goldys, Schoolo of Mathematics and Statistics, University of New
South Wales, Sydney, Australia. Email: \texttt{B.Goldys at
unsw.edu.au}.} \and Fausto Gozzi\footnote{\thinspace Fausto Gozzi
(corresponding author), Dipartimento di Scienze Economiche ed
Aziendali, Facolt\`{a} di Economia, Libera Universit\`{a}
internazionale degli studi sociali ``Guido Carli'', viale Romania
32, 00197 Roma, Italy. Email: \texttt{fgozzi at luiss.it}.} }

\maketitle
\begin{abstract}
This paper, which is the natural continuation of \cite{fegogo},
studies a class of optimal control problems with state constraints
where the state equation is a differential equation with delays.
This class includes some problems arising in economics, in
particular the so-called models with time to build. In \cite{fegogo}
the problem is embedded in a suitable Hilbert space $H$ and the
regularity of the associated Hamilton-Jacobi-Bellman (HJB) equation
is studied. Therein the main result is that the value function $V$
solves the HJB equation and has continuous classical derivative in
the direction of the ``present''. The goal of the present paper is
to exploit such result to find optimal feedback strategies for the
problem. While it is easy to define formally a feedback strategy in
classical sense the proof of its existence and of its optimality is
hard due to lack of full regularity of $V$ and to the infinite
dimension. Finally, we show some approximation results that allow us
to apply our main theorem to obtain $\varepsilon $-optimal
strategies for a wider class of problems.
\end{abstract}


\bigskip

\noindent \textbf{A.M.S. Subject Classification}: 34K35,49L25,
49K25.

\newpage

\tableofcontents

\section{Introduction}

The main purpose of this paper is to prove the existence of optimal
feedback strategies for a class of optimal control problems of
deterministic delay equations arising in economic models.

The paper represents the natural continuation  of \cite{fegogo}
where a class of optimal control problems with state constraints
where the state equation is a differential equation with delays is
studied. This class includes some problems arising in economics, in
particular the so-called models with time to build. In \cite{fegogo}
the problem is embedded in a suitable Hilbert space $H$ and the
associated Hamilton-Jacobi-Bellman (HJB) equation is studied.
Therein the main result is concerned  the regularity of solutions to
such a HJB equation. More precisely it is shown that  the value
function has continuous classical derivative in the direction of the
``present''. This allows to define a feedback strategy in classical
sense.

In the present paper we start from this result and we exploit it to
prove:
\begin{itemize}
    \item the existence of optimal feedback strategies through a Verification
 Theorem;
    \item the existence of $\varepsilon$-optimal strategies for a wider family
of problems through approximation results.
\end{itemize}

The class of optimal control problems is the following: given a
control $c\ge 0$ the state $x$ satisfies the following delay
equation
$$
\begin{cases}
x'(t)=rx(t)+f_0\left(x(t),\int_{-T}^0a(\xi)x(t+\xi)d\xi\right)-c(t),\\
x(0)=\eta_0, \ x(s)=\eta_1(s), \ s\in[-T,0),
\end{cases}
$$
with state constraint $x(\cdot)>0$.
The objective is to maximize the functional
$$J(\eta; c(\cdot)):=\int_0^{+\infty} e^{-\rho t}
\left[U_1(c(t))+ U_2(x(t))\right]\,dt, \ \ \ \rho> 0,$$ over the set
of the admissible controls $c$.

When the feedback strategy effectively exists and is admissible we
prove (Theorem \ref{TH:ver}) that it must be optimal for the
problem: this is not trivial since we do not have the full gradient
of the value function and so we need to use a verification theorem
for viscosity solution which is new in this context. Indeed a
verification theorem in the framework of viscosity solution is given
in the finite dimensional case in \cite{YongZhou}. Adapting the
technique of proof to our case is difficult due to the infinite
dimensional nature of our problem and to a mistake in the key Lemma
5.2, Chapter 5 of \cite{YongZhou} that we discovered here and that
is pointed out in Remark \ref{rem:lemma}. We then give (Proposition
\ref{peano}) sufficient conditions under which the formal optimal
feedback exists and is admissible.

Since our setting (where we prove the Verification Theorem and the
existence of optimal feedback strategies) do not cover the case of
pointwise delay (see \cite{fegogo}, Remark 4.8) which is used in the
previously quoted applications, we go further showing three
approximation results that allow to apply our main theorem to obtain
$\varepsilon $-optimal strategies for a wider class of problems
including the case of pointwise delay (Propositions
\ref{PR:nostateut}, \ref{propV_n,0}, \ref{PR:both}).

\bigskip

The plan of the paper is as follows. Section \ref{Sec:pr} is devoted
to recall the problem set in \cite{fegogo} and the results contained
therein. Section 3 contains the first main result, i.e. the
Verification Theorem  \ref{TH:ver}.  Then, in Section 3.1 we give
sufficient conditions under which the hypothesis of the verification
theorem is satisfied. Section \ref{sec:appr} closes the paper with
the announced approximation results.

\section{The optimal control problem}\label{Sec:pr}
In this section we give the setup of the optimal control problem and
recall, for the reader's convenience, the main results of
\cite{fegogo}. We will use the notations
\[L^2_{-T}:=L^2([-T,0];\mathbb{R}),\quad\mathrm{and}\quad
W^{1,2}_{-T}:=W^{1,2}([-T,0];\mathbb{R}).\]
We will denote by $H$ the Hilbert space
$$H:=\mathbb{R}\times L^2_{-T},$$ endowed with the inner product
$$\langle\cdot,\cdot\rangle= \langle\cdot,\cdot\rangle_{\mathbb{R}}+
\langle\cdot,\cdot\rangle_{L^2_{-T}},$$
and the norm
\[\|\cdot\|^2=|\cdot|^2_{\mathbb R}+\|\cdot\|^2_{L^2_{-T}}.\]
We will denote by
$\eta:=(\eta_0,\eta_1(\cdot))$ the generic element of this space.
For convenience we set also$$H_+:=(0,+\infty)\times L^2_{-T}, \ \ \ \ \
H_{++}:=(0,+\infty)\times\{\eta_1(\cdot)\in L^2_{-T} \ | \
\eta_1(\cdot)\geq 0 \ a.e.\}.$$
\begin{remark}\label{rm:meglioH+}
Economic motivations we are mainly interested in (see
\cite{AseaZak,BambiJEDC,KydlandPrescott})
require to study the optimal control problem with the initial
condition in $H_{++}$. However, the set $H_{++}$ is not convenient
to work with, since its interior with respect to the
$\|\cdot\|$-norm is empty. That is why we enlarge the problem and
allow the initial state belonging to the class $H_+$.
\hfill$\blacksquare$
\end{remark}

For $\eta\in H_+$, we consider an optimal control of the following
differential delay equation:
\begin{equation}\label{eqstate}
\begin{cases}
x'(t)=rx(t)+f_0\left(x(t),\int_{-T}^0a(\xi)x(t+\xi)d\xi\right)-c(t),\\
x(0)=\eta_0, \ x(s)=\eta_1(s), \ s\in[-T,0),
\end{cases}
\end{equation}
with state constraint $x(\cdot)>0$ and control constraint $c(\cdot)\geq 0$.
We set up the following assumptions on the functions $a,f_0$.
\begin{hypothesis}\label{f_0,a}
\begin{itemize}
\item[]
\item $a(\cdot)\in W^{1,2}_{-T}$ is such that $a(\cdot)\geq 0$ and $a(-T)=0$;
\item $f_0:[0,\infty)\times \mathbb{R}\rightarrow \mathbb{R}$
 is jointly concave, {nondecreasing} with respect to the second variable, Lipschitz continuous with Lipschitz constant $C_{f_0}$, and
 \begin{equation}\label{f_0strict}
 f_0(0,y)>0, \ \ \ \ \forall y>0.
 \end{equation}
\end{itemize}
\hfill $\blacksquare$
\end{hypothesis}

\begin{remark}\label{rm:Bambipossibile}
In the papers \cite{AseaZak,BambiJEDC,KydlandPrescott} the pointwise
delay is used. We cannot treat exactly this case for technical
reason that are explained in Remark \ref{rm:noBambi} below. However
we have the freedom of choosing the function $a$ in a wide class and
this allows to take account of various economic phenomena. Moreover
we can approximate the pointwise delay with suitable sequence of
functions $\{a_n\}$ getting convergence of the value function and
constructing $\varepsilon$-optimal strategies (see Subsections
\ref{subs:apprBambi} and \ref{subs:apprBambibis}).
\hfill$\blacksquare$
\end{remark}

\par\medskip\noindent
We say that a {function} $x:[-T,\infty)\to\mathbb R^+$ is a solution to equation \eqref{eqstate} if $x(t)=\eta_1(t)$ for $t\in[-T,0)$ and
\[x(t)=\eta_0+\int_0^trx(s)ds+\int_0^tf_0\left(x(s),\int_{-T}^0a(\xi)x(s+\xi)d\xi\right)ds-\int_0^tc(s)ds,\quad t\ge 0.\]
\begin{theorem}\label{daprato}
 For any given $\eta \in H_+$, $c(\cdot)\in L^1_{loc}([0,+\infty);\mathbb{R^+})$,  equation \eqref{eqstate} admits a unique solution that is absolutely continuous {on $[0,+\infty)$}.
\end{theorem}
\textbf{Proof.} See \cite{fegogo}.\hfill$\square$\\

We denote by $x(\cdot; \eta, c(\cdot))$ the unique solution of
(\ref{eqstate}) with initial point $\eta\in H_+$ and under the
control $c(\cdot)$. We emphasize that this solution actually
satisfies pointwise only the integral equation associated with
\eqref{eqstate}; it satisfies \eqref{eqstate} in differential form
only for almost every $t\in[0,+\infty)$. \vspace{.3cm}

For $\eta\in H_+$ we define the class of the admissible controls
starting from $\eta$ as $$\mathcal{C}(\eta):= \{ c(\cdot)\in
L^1_{loc}([0,+\infty);\mathbb{R^+}) \ | \ x(\cdot;  \eta, c(\cdot))
> 0\}.$$ Setting $x(\cdot):=x(\cdot,\,; \eta,c(\cdot))$, the problem
consists in maximizing the functional
$$J(\eta; c(\cdot)):=\int_0^{+\infty} e^{-\rho t} \left[U_1(c(t))+ U_2(x(t))\right]\,dt, \ \ \ \rho> 0,$$
over the set of the admissible strategies. \vspace{.3cm}

The following will be standing assumptions on the utility functions
$U_1$, $U_2$, holding throughout the whole paper.
\begin{hypothesis}\label{hyp:utility}
\begin{itemize}
\item[]
\item[(i)] $U_1\in C([0,+\infty);\mathbb{R})\cap C^2 ((0,+\infty);\mathbb{R})$,  $U_1'>0$, $U_1'(0^+)=+\infty$, $U_1''<0$ and $U_1$ is bounded.
\item[(ii)] $U_2\in C((0,+\infty);\mathbb{R})$ is increasing, concave, bounded from above. Moreover
\begin{equation}\label{ipotesiU_2}
\int_0^{+\infty} e^{-\rho t} U_2\left(e^{-C_{f_0} t}\right)dt>-\infty.
\end{equation}
\end{itemize}
\vspace{-0,50truecm}\hfill $\blacksquare$
\end{hypothesis}
Since $U_1$, $U_2$ are bounded from above, the previous functional is well-defined for any $\eta\in H_+$ and $c(\cdot)\in\mathcal{C}(\eta)$.
We set
$$\bar{U}_1:=\lim_{s\rightarrow+\infty} U_1(s), \ \ \ \bar{U}_2:=\lim_{s\rightarrow+\infty} U_2(s).$$

We refer to \cite{fegogo} for comments on the assumptions above.


\bigskip
 For $\eta\in H_+$ the value function of our problem is defined by
\begin{equation}\label{value}
V(\eta):=\sup_{c(\cdot)\in\mathcal{C}(\eta)} J(\eta, c(\cdot)),
\end{equation}
with the convention $\sup \emptyset =-\infty.$
 The domain of the value function is the set
$$\mathcal{D}(V):=\{\eta\in H_+ \ | \ V(\eta)>-\infty\}.$$
Due to the assumptions on $U_1$, $U_2$ we directly get that $V\leq \frac{1}{\rho}(\bar{U}_1+\bar{U}_2)$.

\subsection{Preliminary results}\label{sec:value}
The proof of the following qualitative results on the value function
can be found in \cite{fegogo}.
\begin{lemma}[Comparison]\label{comparison}
 Let $\eta\in H_+$ and let $c(\cdot)\in L^1_{loc}([0,+\infty);\mathbb{R^+})$. Let $x(t)$, $t\geq 0$, be an absolutely continuous function satisfying almost everywhere the differential inequality
\begin{equation*}
\begin{cases}
x'(t)\leq rx(t)+f_0\left(x(t),\int_{-T}^0a(\xi)x(t+\xi)d\xi\right)-c(t),\\
x(0)\leq\eta_0, \ x(s)\leq\eta_1(s), \ \mbox{for a.e.} \ s\in[-T,0).
\end{cases}
\end{equation*}
Then
$x(\cdot)\leq x(\cdot;{\eta},{c}(\cdot)).$\hfill$\square$

\end{lemma}

\begin{proposition}\label{piu}
We have
$$H_{++}\subset\mathcal{D}(V), \ \ \ \ \ \mathcal{D}(V)=\{\eta\in H_+ \ | \ 0\in \mathcal{C}(\eta)\}.$$
 The set $\mathcal{D}(V)$ is convex and
the value function $V$  is  concave on $\mathcal{D}(V)$. \hfill$\square$
\end{proposition}

\begin{proposition}\label{strict}
We have the following statements:
\begin{enumerate}
\item $V(\eta)< \frac{1}{\rho}(\bar{U}_1+\bar{U}_2)$ for any $\eta\in H_+$.
\item  $\lim_{\eta_0\rightarrow +\infty} V(\eta_0,\eta_1(\cdot))=\frac{1}{\rho}(\bar{U}_1+\bar{U}_2)$, for all $\eta_1(\cdot)\in L^2_{-T}$.
\item $V$ is strictly increasing with respect to the first variable.
\end{enumerate}
\hfill$\square$
\end{proposition}

\subsection{The delay problem rephrased in  infinite dimension}\label{sec:infinite}
Our aim is to apply the dynamic programming technique in order to solve the control problem described in the previous section. However, this approach requires a markovian setting. That is why we will reformulate the problem as an infinite-dimensional control problem. Let $\hat{n}=(1,0)\in H_+$ and
let us consider, for $\eta\in H$ and $c(\cdot)\in L^1([0,+\infty);\mathbb{R}^+)$, the following evolution equation in the space $H$:
\begin{equation}\label{infinitestate}
\begin{cases}
X'(t)=AX(t)+ F(X(t))-c(t)\hat{n},\\
X(0)=\eta\in H_+.
\end{cases}
\end{equation}
In the equation above:
\begin{itemize}
\item $A:\mathcal{D}(A)\subset H\longrightarrow H$ is an unbounded operator defined by $A(\eta_0, \eta_1(\cdot)):= (r\eta_0, \eta_1'(\cdot))$ on
$$\mathcal{D}(A):=\{ \eta\in H \ | \ \eta_1(\cdot)\in W^{1,2}_{-T}, \ \eta_1(0)=\eta_0\};$$
\item $F:H\longrightarrow H$ is a Lipschitz continuous map defined by
\begin{equation*}
F(\eta_0,\eta_1(\cdot)):=\left(f\left(\eta_0,\eta_1(\cdot)\right),0\right),
\end{equation*}
where $f(\eta_0,\eta_1(\cdot)):=f_0\left(\eta_0,\int_{-T}^0a(\xi)\eta_1(\xi)d\xi\right)$.
\end{itemize}
It is well known that $A$ is the infinitesimal generator of a
strongly continuous semigroup $(S(t))_{t\geq 0}$ on $H$; its explicit expression is given by
$$
S(t)(\eta_0,\eta_1(\cdot))=\left(\eta_0e^{rt},I_{[-T,0]}(t+\cdot) \
\eta_1(t+\cdot)+ I_{[0,+\infty)}(t+\cdot) \ \eta_0e^{r(t+\cdot)}\right);
$$

\subsubsection{Mild solutions of the state equation}\label{subs:mild}
Here we give a definition of the \emph{mild solution}
to  (\ref{infinitestate}),  state the existence and uniqueness of
such a solution and the equivalence between the one dimensional
delay problem and the infinite dimensional one. We refer to \cite{fegogo} for the proofs.
\begin{definition}
A mild solution of
(\ref{infinitestate})
is a function $X\in C([0,+\infty);H)$ which satisfies the integral equation
\begin{equation}\label{mildcap3}
X(t)=S(t)\eta+\int_0^t S(t-\tau)F(X(\tau))d\tau
+\int_0^t c(\tau)S(t-\tau)\hat{n}\,d\tau.
\end{equation}
\end{definition}
\begin{theorem}\label{existmild}
For any $\eta\in H$, there exists a unique mild solution of (\ref{infinitestate}). \hfill$\square$
\end{theorem}

We denote by $X(\cdot;\eta, c(\cdot))=(X_0(\cdot;\eta,c(\cdot)), X_1(\cdot;\eta,c(\cdot)))$ the unique solution to (\ref{infinitestate}) for the initial state $\eta\in H$ and under the control $c(\cdot)\in L^1([0,+\infty);\mathbb{R}^+)$.
The following equivalence result justifies our approach.
\begin{proposition} \label{link}
Let $\eta\in H_+$, $c(\cdot)\in\mathcal{C}(\eta)$ and let $x(\cdot)$, $X(\cdot)$ be respectively the unique solution to (\ref{eqstate}) and the unique mild solution to (\ref{infinitestate}) starting from $\eta$ and under the control $c(\cdot)$. Then, for any $t\geq 0$, we have the equality in $H$  $$X(t)=\left(x(t),x(t+\xi)_{\xi\in [-T,0]}\right).$$ \hfill$\square$
\end{proposition}

\subsubsection{Regularity of the value function}\label{sec:continuity}
Here  we state the regularity properties  of the value function. We refer to \cite{fegogo} for the proofs.

We recall that the generator $A$ of the semigroup $(S(t))_{t\geq 0}$ has bounded inverse in $H$ given by
\[A^{-1}\left(\eta_0,\eta_1\right)(s)=\left(\frac{\eta_0}{r},\frac{\eta_0}{r}-\int_s^0\eta_1(\xi)d\xi\right),\quad s\in [-T,0].\]
It is well known that $A^{-1}$ is compact in $H$. It is also clear that $A^{-1}$ is an isomorphism of $H$ onto $\mathcal{D}(A)$ endowed with the graph norm.

We define the $\|\cdot\|_{-1}$-norm on $H$ by $$\|\eta\|_{-1}:=\|A^{-1}\eta\|.$$

\begin{proposition}\label{usefcap3}
The set $\mathcal{D}(V)$ is open in the space $\left(H,\|\cdot\|_{-1}\right)$ and the value function is continuous with respect to $\|\cdot\|_{-1}$ on $\mathcal{D}(V)$. Moreover
\begin{equation}\label{debcont}
(\eta_n)\subset\mathcal{D}(V), \ \  \eta_n\rightharpoonup \eta\in\mathcal{D}(V) \ \ \Longrightarrow \ \ V(\eta_n)\rightarrow V(\eta).
\end{equation}\hfill$\square$
\end{proposition}

Therefore, we can apply the following result to the value function.
\begin{proposition}\label{importante}
 Let $v:\mathcal{D}(V)\rightarrow \mathbb{R}$ be a concave function continuous with respect to $\|\cdot\|_{-1}$. Then
\begin{enumerate}
\item $v=u\circ A^{-1}$, where
$u:\mathcal{O}\subset H \rightarrow \mathbb{R}$
is a concave $\|\cdot\|$-continuous function.
\item $D^+v(\eta)\subset \mathcal{D}(A^*)$, for any $\eta\in \mathcal{D}(V)$.
\item $D^+ u(A^{-1}\eta)= A^*D^+v(\eta)$, for any $\eta\in\mathcal{D}(V)$. In particular, since $A^*$ is injective, $v$ is differentiable at $\eta$ if and only if $u$ is differentiable at $A^{-1}\eta$.
\item If $\zeta\in D^*v(\eta)$, then there exists a sequence $\eta_n\rightarrow \eta$ such that there exist $\nabla v(\eta_n)$, $\nabla v(\eta_n)\rightarrow \zeta$ and $A^*\nabla v(\eta_n)\rightharpoonup A^*\zeta$.
\end{enumerate}
\hfill$\square$
\end{proposition}

The HJB equation associated to our optimization problem is
\begin{equation}\label{HJB}
\rho v(\eta)= \langle \eta, A^* \nabla v(\eta)\rangle +f(\eta)v_{\eta_0}(\eta)+ U_2(\eta_0)+\mathcal{H}(v_{\eta_0}(\eta)),
\end{equation}
where $\mathcal{H}$ is the Legendre transform of $U_1$, i.e.
$$\mathcal{H}(\zeta_0):=\sup_{c\geq 0}\left(U_1(c)-\zeta_0 c\right), \  \ \ \zeta_0> 0.$$
Due to  Hyphothesis \ref{hyp:utility}-(i) and to Corollary 26.4.1 of \cite{Rock}, we have that $\mathcal{H}$ is strictly convex on $(0,+\infty)$. Notice that, thanks to Proposition \ref{strict}-(3),
$$D^+_{\eta_0}V(\eta):=\{\zeta_0 \in \mathbb{R} \ | \ (\zeta_0,\zeta_1(\cdot))\in D^+V(\eta)\}\subset (0,\infty)$$
for any $\eta\in\mathcal{D}(V)$, i.e. where $\mathcal{H}$ is defined.

We can study this equation following the viscosity approach. In order to do that, we have to define a
suitable set of regular test functions. This is the set
\begin{equation}\label{testcap3}
\tau:=\Big\{\varphi\in C^1(H) \ | \ \nabla \varphi (\cdot)\in \mathcal{D}(A^*), \ \eta_n\rightarrow \eta\Rightarrow A^*\nabla\varphi(\eta_n)\rightharpoonup A^*\nabla\varphi(\eta)\Big\}.
\end{equation}
Let us define, for $c\geq 0$,  the operator $\mathcal{L}^{c}$ on $\tau$ by
$$[\mathcal{L}^{c}\varphi](\eta):=-\rho \varphi (\eta)+\langle \eta ,A^*\nabla\varphi(\eta)\rangle+  f(\eta)\varphi_{\eta_0}(\eta)-c\varphi_{\eta_0}(\eta).$$

\begin{definition}\label{def:viscosity}

 (i) A continuous function $v:\mathcal{D}(V)\rightarrow \mathbb{R}$ is called a viscosity subsolution of (\ref{HJB}) on $\mathcal{D}(V)$ if for any $\varphi\in\tau$ and any $\eta_M\in\mathcal{D}(V)$ such that $v-\varphi$ has a $\|\cdot\|$-local maximum at $\eta_M$ we have
$$\rho v(\eta_M)\leq   \langle \eta_M, A^* \nabla \varphi(\eta_M)\rangle +f(\eta_M)\varphi_{\eta_0}(\eta_M)+U_2(\eta_0)+ \mathcal{H}(\varphi_{\eta_0}(\eta_M)).$$

(ii) A continuous function $v:\mathcal{D}(V)\rightarrow \mathbb{R}$ is called a viscosity supersolution of (\ref{HJB}) on $\mathcal{D}(V)$ if for any $\varphi\in\tau$ and any $\eta_m\in\mathcal{D}(V)$ such that $v-\varphi$ has a $\|\cdot\|$-local minimum at $\eta_m$  we have
$$\rho v(\eta_m)\geq   \langle \eta_m, A^* \nabla \varphi(\eta_m)\rangle +f(\eta_m)\varphi_{\eta_0}(\eta_m)+ U_2(\eta_0)+ \mathcal{H}(\varphi_{\eta_0}(\eta_m)).$$

(iii) A continuous function $v:\mathcal{D}(V)\rightarrow \mathbb{R}$ is called a viscosity supersolution of (\ref{HJB}) on $\mathcal{D}(V)$ if it is both a viscosity sub and supersolution.

\end{definition}

\begin{theorem}\label{TH:visc}
The value function $V$ is a viscosity solution of (\ref{HJB}) on $\mathcal{D}(V)$. \hfill$\square$
\end{theorem}

Actually  the concave $\|\cdot\|_{-1}$-continuous viscosity solutions of (\ref{HJB}) (so that in particular the value function $V$) are differentiable along the direction $\hat{n}=(1,0)$. This is stated in the next result: we refer to \cite{fegogo} for the proof.
\begin{theorem}\label{TH:reg}
 Let $v$ be a concave $\|\cdot\|_{-1}$-continuous viscosity solution of (\ref{HJB}) on $\mathcal{D}(V)$. Then $v$ is differentiable along the direction $\hat{n}=(1,0)$ at any point $\eta\in\mathcal{D}(V)$ and the function $\eta\mapsto v_{\eta_0}(\eta)$ is continuous on $\mathcal{D}(V)$. \hfill $\square$
\end{theorem}

\section{Verification theorem and optimal feedback strategies}\label{secopt}

Here we prove a Verification Theorem yielding optimal strategies for the problem. We start with the following definition.
\begin{definition}
Let $\eta\in \mathcal{D}(V)$. An admissible control $c^*(\cdot)\in\mathcal{C}(\eta)$ is said to be optimal for the initial state $\eta$ if $J(\eta; c^*(\cdot))=V(\eta)$. In this case the corresponding state trajectory $x^*(\cdot):=x(\cdot; \eta, c^*(\cdot))$ is said to be an optimal trajectory and the couple $(x^*(\cdot),c^*(\cdot))$ is said an optimal couple.
   \hfill$\square$
\end{definition}

Thanks to the regularity result of the previous section 
 we can
define, at least formally, the ``candidate'' optimal feedback map on
$\mathcal{D}(V)$, which is given by
\begin{equation}\label{feedbackmap}
C(\eta):=\mbox{argmax}_{c\geq 0} \left(U_1(c)-c V_{\eta_0}(\eta)\right), \ \ \ \eta\in \mathcal{D}(V).
\end{equation}
Note that
this map  is well-defined 
since $V$ is concave and, by Proposition
\ref{strict},  strictly increasing,  so that  we have $V_{\eta_0}(\eta)\in
(0,+\infty)$ for all $\eta\in\mathcal{D}(V)$. Existence and
uniqueness of the argmax follow from the assumptions on $U_1$.
Moreover, since $V_{\eta_0}$ is continuous on $\mathcal{D}(V)$, also
$C$ is continuous on $\mathcal{D}(V)$. The closed-loop delay state
equation associated with this map is, for $\eta\in\mathcal{D}(V)$,
\begin{equation}\label{rfe}
\begin{cases}
x'(t)=rx(t)+f_0\left(x(t),\int_{-T}^0a(\xi)x(t+\xi)d\xi\right)-
{C}\left((x(t),x(t+\xi)|_{\xi\in[-T,0]})\right),\\
x(0)=\eta_0, \ x(s)=\eta_1(s), \ s\in[-T,0).
\end{cases}
\end{equation}
Now we want to prove a Verification Theorem: if the closed loop
equation \eqref{rfe} has a strictly positive solution $x^*(\cdot )$,
(so that we must have
$\left(x^*(t),x^*(t+\xi)|_{\xi\in[-T,0]}\right)\in\mathcal{D}(V)$ and the term
$C\left((x^*(t),x^*(t+\xi)|_{\xi\in[-T,0]})\right)$ is well-defined
for every $t\geq 0$),
then the feedback strategy
\begin{equation}\label{stratopt} c^*(t):=
C\left((x^*(t),x^*(t+\xi)|_{\xi\in[-T,0]})\right)
\end{equation}
is optimal.
Notice that, by definition of $c^*(\cdot)$, if  $x^*(\cdot)$ is a strictly positive solution of \eqref{rfe}, then $c^*(\cdot)$ is admissible and, setting $X^*(t):=X(t;\eta,c^*(\cdot))$, we have
$$X^*(t)=\left(x^*(t),x^*(t+\xi)|_{\xi\in[-T,0]}\right)\in\mathcal{D}(V), \ \ \ \ \forall t\geq 0.$$
In order to prove a Verification Theorem, formally we need to integrate the function
\begin{equation}\label{optraj}
t\mapsto \frac{d}{dt}\left[e^{-\rho t}V(X^*(t))\right].
\end{equation}
Thus we need something like a Fundamental Theorem of Calculus relating the function and the integral of its "derivative". Since we do not require the intial datum $\eta$ belonging to $\mathcal{D}(A)$ and the operator $A$ works as a shift operator on the infinite-dimensional component, we do not have the condition $X^*(t)\in\mathcal{D}(A)$ for almost every $t\geq 0$ giving a regularity for the function
$$t\mapsto e^{-\rho t} V(X^*(t))$$
sufficient to apply the Fundamental Theorem of Calculus (see \cite{LiYong}, Theorems 5.4, 5.5, Chapter 6). Therefore we can suppose only that the function \eqref{optraj} is continuous and we should try to apply a generalized Fundamental Theorem of Calculus in inequality form. There is such a result in \cite{YongZhou}, Lemma 5.2, Chapter 5. Unfortunately such result is not true as it is stated (see Remark \ref{rem:lemma} for a counterexample), so we have to refer to other results based on the theory first formulated by Dini and Lebesgue. We refer to \cite{Bru, HagTho, Saks} sketching the ideas we need. If $g$ is a continuous function on some interval $[\alpha,\beta]\subset\mathbb{R}$, the right Dini derivatives of $g$   are defined by
$$D^+g(t)=\limsup_{h\downarrow 0} \frac{g(t+h)-g(t)}{h}, \ \ D_+g(t)=\liminf_{h\downarrow 0} \frac{g(t+h)-g(t)}{h}, \ \ \ t\in[\alpha, \beta),$$
and the left Dini derivatives by
$$D^-g(t)=\limsup_{h\uparrow 0} \frac{g(t+h)-g(t)}{h}, \ \ D_-g(t)=\liminf_{h\uparrow 0} \frac{g(t+h)-g(t)}{h}, \ \ \ t\in(\alpha,\beta].$$
The key result is the following (see \cite{Bru}, Theorem 1.2, Chapter 4).
\begin{proposition}\label{prop:bru}
If $g$ is a continuous real function on $[\alpha,\beta]$, then the bounds of  each  Dini's derivative are equal to the bounds of the set of the difference quotients
$$\left\{\frac{g(t)-g(s)}{t-s} \  \Big| \ t,s\in[\alpha,\beta]\right\}.$$
\hfill$\square$
\end{proposition}
\noindent
An immediate consequence of Proposition \ref{prop:bru} above  is the following.
\begin{proposition}[Monotonicity result]\label{monotonicityresult}
Let $g\in C([\alpha,\beta];\mathbb{R})$ be such that $D^+g(t)\geq 0$ for all $t\in[\alpha, \beta)$. Then $g$ is nondecreasing on $[\alpha, \beta]$.\hfill $\square$
\end{proposition}

Now we can give a simple lemma useful for proving the Verification Theorem.\begin{lemma}\label{lemmagiorgio}
Let $g,\mu\in C([0,+\infty);\mathbb{R})$ 
such that
\begin{equation}\label{everywhere}
D_-g(t)\geq \mu(t), \ \ \  \forall t\in(0,+\infty).
\end{equation}
Then, for every $0\leq \alpha\leq\beta<+\infty$,
\begin{equation}\label{tesilemma}
g(\beta)-g(\alpha)\geq \int_{\alpha}^{\beta} \mu(t) dt.
\end{equation}
\hfill$\square$
\end{lemma}
\textbf{Proof.} Since $D_-g(t)\geq \mu(t)$ for every
$t\in(0,+\infty)$, we have $D_-[g(t)-\int_0^t \mu(s)ds]\geq 0$ for
every $t\in(0,+\infty)$. Thanks to Proposition \ref{prop:bru} we
have also $D^+[g(t)-\int_0^t \mu'(s)ds]\geq 0$ for every
$t\in[0,+\infty)$. Therefore, due to Proposition
\ref{monotonicityresult}, $t\mapsto g(t)-\int_0^t\mu (s)ds$ is
nondecreasing, getting the claim.\hfill $\square$

\begin{remark}\label{rem:lemma}
Following \cite{HagTho}, we give some remarks on Lemma \ref{lemmagiorgio}.
\begin{itemize}
\item
The assumption that $\mu$ is continuous can be replaced assuming
that  $\mu$ is a finite-valued (Lebesgue) measurable and integrable function (Theorem 9
of \cite{HagTho}); also  condition \eqref{everywhere} can be weakened  assuming that it holds
out of a countable set (Section 5.b of \cite{HagTho}).
\item
Condition \eqref{everywhere} can be weakened assuming that it holds
almost everywhere adding the assumption $D_-g>-\infty$ everywhere
(Section 5.c of \cite{HagTho}).
\item
We cannot further weaken \eqref{everywhere}: if  it is verified only
almost everywhere without any further assumption on $D_-g$, then
\eqref{tesilemma} is no longer true. For example, if $g=-f$ on
$[0,1]$, where $f$ is the Cantor function and $\mu\equiv 0$, we have
$$\mu(t)=0=g'(t)=D_-g(t) \ \ \ \mbox{for a.e.} \ t\in(0,1].$$
Therefore,  taking $\alpha=0, \ \beta=1$, the left handside of
\eqref{tesilemma} is $-1$, while the right handside is $0$. Indeed
in this case $D_-g=-\infty$ on the Cantor set. So Lemma 5.2, Chapter
5, of \cite{YongZhou} is not correct. Indeed the condition required
therein is not sufficient to apply Fatou's Lemma in the proof.
\hfill$\blacksquare$
\end{itemize}
\end{remark}
\begin{theorem}[Verification]\label{TH:ver}
Let $\eta\in H_+$ and let $x^*(\cdot)$ be a solution of \eqref{rfe} such that $x^*(\cdot)>0$; let $c^*(\cdot)$ be the strategy defined by \eqref{stratopt}. Then $c^*(\cdot)$ is admissible and optimal for the problem.
\end{theorem}
\textbf{Proof.} As said above the fact that $c^*(\cdot)$ is
admissible is a direct consequence of the assumption $x^*(\cdot)>0$
and of the definition of $c^*(\cdot)$.

Set $X^*(\cdot):=X(\cdot\,;\eta,c^*(\cdot))$ and let $s>0$. Let
$p_1(s)\in L^2_{-T}$ be such that $$\left(V_{\eta_0}(X^*(s)),
p_1(s)\right)\in D^+ V(X^*(s))$$ and let
$$\varphi(\zeta):=V(X^*(s))+\langle  \left(V_{\eta_0}(X^*(s)), p_1(s)\right), \zeta- X^*(s)\rangle, \ \ \ \ \zeta\in H,$$
so that
$$\varphi (X^*(s))=V(X^*(s)), \ \ \ \ \ \varphi (\zeta)\geq V(\zeta), \ \ \zeta\in H.$$
From Proposition \ref{importante} we know that $\varphi\in\tau$, so that
\begin{multline*}
\liminf_{h\uparrow 0}\frac{e^{-\rho(s+h)}V(X^*(s+h))-e^{-\rho s}V(X^*(s))}{h}\geq \liminf_{h\uparrow 0}\frac{e^{-\rho(s+h)}\varphi(X^*(s+h))-e^{-\rho s}\varphi(X^*(s))}{h}\\
=e^{-\rho s}\left[\mathcal{L}^{c^*(s)}\varphi\right](X^*(s))
= e^{-\rho s}\Big[-\rho V(X^*(s))+\langle X(s) ,A^* \left(V_{\eta_0}(X^*(s)), p_1(s)\right)\rangle\\+  f(X^*(s))V_{\eta_0}(X^*(s))+c^*(s) V_{\eta_0}(X^*(s))\Big].
\end{multline*}
Due to the definition of $c^*(\cdot)$ we get
\begin{multline*}\liminf_{h\uparrow 0}\frac{e^{-\rho(s+h)}V(X^*(s+h))-e^{-\rho s}V(X^*(s))}{h}+ e^{-\rho s}[U_1(c^*(s))+U_2(X^*_0(s))]\\
\geq   e^{-\rho s}\Big[-\rho V(X^*(s))+\langle X^*(s) ,A^* \left(V_{\eta_0}(X^*(s)), p_1(s)\right)\rangle \\+ f(X^*(s))V_{\eta_0}(X^*(s))+\mathcal{H}(X^*(s))+U_2(X^*_0(s))\Big].
\end{multline*}
Due to the subsolution property of $V$ we get
\begin{equation*}
\label{liminf}
\liminf_{h\uparrow 0}\frac{e^{-\rho(s+h)}V(X^*(s+h))-e^{-\rho s}V(X^*(s))}{h}+ e^{-\rho s}[U_1(c^*(s))+U_2(X^*_0(s))]\geq 0.
\end{equation*}
The function $s\mapsto e^{-\rho s}V(X^*(s))$ and the function $s\mapsto e^{-\rho s}[U_1(c^*(s))+U_2(X^*_0(s))]$ are continuous; therefore we can  apply Lemma \ref{lemmagiorgio}  on  $[0,M]$, $M>0$, getting
$$
e^{-\rho M} V(X^*(M))+\int_0^M e^{-\rho s}[U_1(c^*(s))+U_2(X^*_0(s))]ds\geq V(\eta).
$$
Since $V, U_1, U_2$ are bounded from above, taking the limsup for $M\rightarrow +\infty$ we get
by Fatou's Lemma
$$
\int_0^{+\infty} e^{-\rho s}[U_1(c^*(s))+U_2(X^*_0(s))]ds\geq V(\eta),
$$
which gives the claim.
\hfill$\square$

\begin{remark}
We have given in Theorem \ref{TH:ver} a sufficient condition of
optimality: indeed, we have proved that if the feedback map defines
an admissible strategy then such a  strategy is optimal. Of course,
a natural question arising is whether, at least with a special
choice of data, such a condition is also necessary for the
optimality, i.e. if, given an optimal strategy, it can be written as
feedback of the associated optimal state. From the viscosity point
of view the answer to this question relies in requiring that the
value function is a \emph{bilateral} viscosity subsolution of
\eqref{HJB} along the optimal state trajectory, i.e. requiring that
the value function satisfies the property of Definition
\ref{def:viscosity}-(i) also with the reverted inequality along this
trajectory.

Such a property of the value function is related to the so-called
\emph{backward dynamic programming principle} which is, in turn,
related to the backward study of the state equation (see \cite{BCD},
Chapter III, Section 2.3). Differently from the finite-dimensional
case, this topic  is not standard in infinite-dimension unless the
operator $A$ is the generator of a strongly continuous \emph{group},
which is not our case.

However, in our case we can use the delay original setting of the
state equation to approach this topic. Then the problem reduces to
find, at least for sufficient regular data, a \emph{backward
continuation} of the solution. This problem is faced, e.g., in
\cite{Hale}, Chapter 2, Section 5. Unfortunately our equation does
not fit the main assumption required therein, which  in our setting
basically corresponds to require that the function $a(\cdot)$, seen
as measure, has an atom at $-T$. Investigation on this is left for
future research. \hfill$\blacksquare$
\end{remark}
\smallskip

\subsection{The closed loop equation}

Up to now we did not make any further assumption on the functions
$a$ and  $U_2$ beyond Hypotheses \ref{f_0,a} and \ref{hyp:utility};
in particular it could be $U_2\equiv 0$. However without any further
assumption  we have no information on the behaviour of $V_{\eta_0}$
when we approach the boundary of $\mathcal{D}(V)$  and therefore we
are not able to say anything about the existence of solutions of the
closed loop equation and whether they satisfy or not the state
constraint. So basically we cannot say whether the hypothesis  of
Theorem \ref{TH:ver} is satisfied or not. In order to give
sufficient conditions for that, we need to do some further
assumptions.

\begin{hypothesis}
We will make use of the following assumptions
\begin{equation}\label{a>0}
(i) \ U_2 \ \mbox{is not integrable at} \ 0^+, \ \ \ (ii) \
\int_{-\varepsilon}^0a(\xi)d\xi>0, \ \ \ \forall\varepsilon>0.
\end{equation}
\hfill$\blacksquare$
\end{hypothesis}

Also we need the following Lemma; we refer to \cite{fegogo} for the proof.
\begin{lemma}\label{estimateA}
Let $X(\cdot), \bar{X}(\cdot)$ be the mild solutions to (\ref{infinitestate}) starting respectively from $\eta, \bar{\eta}\in H$ and both under the null control. Then there exists a constant $C>0$ such that
$$\|X(t)-\bar{X}(t)\|_{-1}\leq C\|\eta-\bar{\eta}\|_{-1}, \ \ \forall t\in[0,T].$$
In particular
$$|X_0(t)-\bar{X}_0(t)|\leq rC\|\eta-\bar{\eta}\|_{-1}, \ \ \forall t\in[0,T].$$
\hfill$\square$
\end{lemma}



\begin{lemma}\label{lemma:boundarybehaviour}
\begin{enumerate}
\item[]
\item The following holds
$$
\partial_{\|\cdot\|} \mathcal{D}(V)=\partial_{\|\cdot\|_{-1}} \mathcal{D}(V).
$$
Thanks to the previous equality we  write without ambiguity  $\partial \mathcal{D}(V)$ for denoting the boundary of $\mathcal{D}(V)$ referred to $\|\cdot\|$ or $\|\cdot\|_{-1}$ indifferentely.
\item
Suppose that \eqref{a>0}-(i) holds; then
$$\lim_{\eta\rightarrow \bar{\eta}} V_{\eta_0}(\eta)=+\infty, \ \ \ \forall\bar{\eta}\in \partial \mathcal{D}(V),$$
where the limit  is taken  with respect to $\|\cdot\|$.
\end{enumerate}
\end{lemma}

\textbf{Proof.} We work with the original one-dimensional state
equation with delay.

\smallskip
1. First of all note that, thanks to Proposition \ref{strict} and
Proposition \ref{usefcap3}, the set $\mathcal{D}(V)$ has the
following structure
\begin{equation}\label{structureofD(V)}
\mathcal{D}(V)=\bigcup_{\eta_1\in L^2_{-T}}\big((\eta_0^{\eta_1},+\infty)\times \{\eta_1\}\big),
\end{equation}
where, for $\eta_1\in L^2_{-T}$, we set $\eta_0^{\eta_1}=\inf
\{\eta_0>0 \ | \ (\eta_0,\eta_1(\cdot))\in \mathcal{D}(V)\}.$
For any $\eta\in H$ set $x^{\eta}(\cdot):=x(\cdot;\eta,0)$
and consider the function $g:H\rightarrow\mathbb{R}$ defined by
$$g(\eta_0,\eta_1(\cdot)):=\inf_{t\in[0,T]} x^{\eta}(t).$$
Thanks to Lemma \ref{estimateA} this function is continuous (with
respect to both the norms $\|\cdot\|$ and $\|\cdot\|_{-1}$), so we
have the following representation of $\mathcal{D}(V)$ in terms of
$g$:
$$\mathcal{D}(V)=\{g>0\}.$$
Lemma  \ref{comparison} shows that $g$ is increasing with respect to
the first variable. Actually $g$ is strictly increasing with respect
to the first variable. Let us show this fact. Let $\eta_1\in
L^2_{-T}$ and take $\eta_0,\bar{\eta}_0\in \mathbb{R}$ such that
$\eta_0>\bar{\eta}_0$. Define
$y(\cdot):=x(\cdot;(\eta_0,\eta_1(\cdot)),0)$,
$x(\cdot):=x(\cdot;(\bar{\eta}_0,\eta_1(\cdot)),0)$ and let
$z(\cdot),\bar{z}(\cdot)$ be respectively the solutions on $[0,T]$
of the  differential problems \emph{without delay}
$$\begin{cases}
z'(t)=rz(t)+f_0\left(z(t),\int_{-T}^0a(\xi)x(t+\xi)d\xi\right),\\
z(0)=\eta_0,
\end{cases}
$$
$$\begin{cases}
\bar{z}'(t)=r\bar{z}(t)+f_0\left(\bar{z}(t),\int_{-T}^0a(\xi)x(t+\xi)d\xi\right),\\
\bar{z}(0)=\bar{\eta}_0,
\end{cases}
$$
Then we have, on the interval $[0,T]$,  $\bar{z}(\cdot)\equiv
x(\cdot)$ and, by comparison criterion, $y(\cdot)\geq z(\cdot)$;
moreover we can apply the classic Cauchy-Lipschitz Theorem for ODEs
getting uniqueness for the solutions of the above ODEs, which yields
$z(\cdot)>\bar{z}(\cdot)$ on $[0,T]$. Thus  $y(\cdot)>x(\cdot)$ on
$[0,T]$, proving that $g$ is strictly increasing with rescpect to
the first variable.

The continuity (with respect to both the norms) of $g$,
\eqref{structureofD(V)} and the fact that $g$ is strictly increasing
with respect to the first variable lead to have
$$
 \partial_{\|\cdot\|} \mathcal{D}(V)= \partial_{\|\cdot\|_{-1}} \mathcal{D}(V)=\{g=0\}=\bigcup_{\eta_1\in L^2_{-T}} \big(\{\eta_0^{\eta_1}\}\times\{\eta_1\}\big).
$$

\smallskip
2. We will intend the topological notions referred  to $\|\cdot\|$. Firstly we prove that
$$\lim_{\eta\rightarrow \bar{\eta}} V(\eta)=-\infty, \ \ \ \forall\bar{\eta}\in \partial \mathcal{D}(V).$$
Let $\bar{\eta}\in \partial\mathcal{D}(V)$ and let $(\eta^n)\subset \mathcal{D}(V)$ be a sequence  such that $\eta^n\rightarrow \bar{\eta}$. We can suppose without loss of generality that $(\eta^n)\subset B(\bar{\eta},1)$.
Set  $$x^n(\cdot):=x(\cdot;\eta^n,0),
 \ \ \ \ p^{n}:=\sup_{\xi\in [0,2T]} x^{n}(\xi).$$
Thanks to Lemma \ref{estimateA} there exists $K>0$ such that $p^n\leq K$ for any $n\in\mathbb{N}$.
So, {since $f_0(x,y)\leq C_0(1+|x|+|y|)$ for some $C_0>0$,  we have for the dynamics of $x^{n}(\cdot)$ in the interval $[0,2T]$
$$\frac{d}{dt}\,x^{n}(t)\leq r x^{n}(t) +R,$$
where
$$R:=C_0    \left(1+K+ \|a\|_{L^2_{-T}}(\|\bar{\eta}_1    \|_{L^2_{-T}}+1)+\|a\|_{L^2_{-T}}T^{1/2}K\right).$$}
Therefore there exists $C>0$ such that, for any $s\in[0,T)$, $n\in\mathbb{N}$,
\begin{equation}\label{ss}
x^{n}(t)\leq x^{n}(s)\, e^{r(t-s)}+\frac{R}{r}\,(e^{r(t-s)}-1)\leq x^n(s) (1+C(t-s))+C(t-s), \ \ \ t\in[s,2T].
\end{equation}
By continuity of $g$ we have $\lim_{n\rightarrow
\infty}g(\eta_0^n,\eta_1^n(\cdot))=0$. Thus for any $\varepsilon>0$
we can find $n_0\in\mathbb{N}$ such that, for $n\geq n_0$, there
exists $s_n\in[0,T)$ such that
\begin{equation}\label{s}
x^{n}(s_n)\leq \varepsilon.
\end{equation}
We want to show that
\begin{equation}\label{mnb}
\int_0^{+\infty}e^{-\rho t} U_2(x^n(t))dt \longrightarrow -\infty, \
\ \ \ n\rightarrow\infty.
\end{equation}
For this purpose, since $U_2$ is bounded from above, it is clear
that we can assume without loss of generality $U_2(\cdot)\leq 0$. We
have for $n\geq n_0$, taking into account \eqref{ss} and \eqref{s},
\begin{eqnarray*}
\int_0^{+\infty}e^{-\rho t} U_2(x^n(t))dt&\leq &
e^{-2\rho T} \int_{s_n}^{2T}U_2\big(x^n(s_n) (1+C(t-s_n))+C(t-s_n)\big)dt\\
&\leq & e^{-2\rho T} \int_{s_n}^{2T}U_2(\varepsilon (1+C(t-s_n))+C(t-s_n))dt\\
&\leq & \frac{e^{-2\rho T}}{C(\varepsilon+1)} \int_{\varepsilon}^{CT}U_2(x)dt.
\end{eqnarray*}
Therefore, by the arbitrariness of $\varepsilon$ and since $U_2$ is
not integrable at $0^+$, we get \eqref{mnb}. This is enough to
conclude that $J(\eta^n;0)\rightarrow -\infty $, as $n\rightarrow
\infty$. Of course we have $x^n(\cdot)\geq x(\cdot;\eta^n,c(\cdot))$
for any $c(\cdot)\in\mathcal {C}(\eta^n)$. Since $U_1$ is bounded
from above this is enough to say that also $V(\eta^n)\rightarrow
-\infty$, as $n\rightarrow \infty$.

Now  we prove the claim.
Let $\bar{\eta}\in \partial\mathcal{D}(V)$ and $(\eta^n)\subset
\mathcal{D}(V)$ be such that $\eta^n\rightarrow \bar{\eta}$, suppose
without loss of generality $(\eta_n)\subset B(\bar{\eta},1)$ and set
$x^n(\cdot):=x(\cdot\,; (\eta^n_0+1, \eta_1^n),0)>0$. Since $f_0$ is
Lipschitz continuous and nondecreasing on the second variable, there
exists $C>0$ such that
$$f_0\left(x(t),\int_{-T}^0 a(\xi)x(t+\xi)d\xi\right)\geq -
C\left(1+x(t)+\|a\|_{L^2_{-T}}(\|\bar{\eta}_1\|_{L^2_{-T}}+1)\right)=:-\tilde{R}.$$
Suppose $\tilde{R}\leq 0$ . Then $\frac{d}{dt} x^n(t)\geq rx^n(t)$,
so that, since $\eta_0^n>0$, we have $x^n(t)\geq \eta_0^n+1\geq 1.$
This leads to the estimate
$$V(\eta_0^n+1,\eta_1^n(\cdot))\geq K, \ \ n\in\mathbb{N},$$
for some $K>0$. Due to the  concavity of $V$ we have the estimate
$$V_{\eta_0}(\eta^n)\geq V(\eta_0^n+1,\eta_1^n(\cdot))- V(\eta_0^n,\eta_1^n(\cdot))\geq K- V(\eta_0^n,\eta_1^n(\cdot))\rightarrow +\infty,$$
i.e. the claim.\\
Suppose then $\tilde{R}>0$ and set $x^n(\cdot):= x(\cdot\,;
(\eta^n_0+\tilde{R}/r, \eta_1^n),0).$ Then $\frac{d}{dt} x^n(t)\geq
rx^n(t)-\tilde{R}$, so that, since $\eta_0^n>0$, we have $x^n(t)\geq
\tilde{R}/r>0$. This leads to the estimate
$$V(\eta_0^n+\tilde{R}/r,\eta_1^n(\cdot))\geq K, \ \ n\in\mathbb{N},$$
for some $K>0$. Due to the  concavity of $V$ we have the estimate
$$V_{\eta_0}(\eta^n)\geq \frac{r}{\tilde{R}}\left[V(\eta_0^n+\tilde{R}/r,\eta_1^n(\cdot))- V(\eta_0^n,\eta_1^n(\cdot))\right]\geq \frac{r}{\tilde{R}}\left[K- V(\eta_0^n,\eta_1^n(\cdot))\right]\rightarrow +\infty,$$
i.e. the claim.
\hfill$\square$\\

\begin{proposition}\label{peano}
Let  \eqref{a>0} hold,
let $\eta\in H_{++}$ and consider the closed-loop delay state equation
\eqref{rfe}.
Then this equation admits a solution $x^*(\cdot)\in C^1([0,+\infty);\mathbb{R})$.
Moreover, for all $t\geq 0$,
$$x^*(t)>0, \ \ \  (x^*(t),x^*(t+\xi)|_{\xi\in[-T,0]})\in\mathcal{D}(V).$$
In particular the feedback strategy defined in
\eqref{stratopt}
is admissible.
\end{proposition}

\textbf{Proof.} Thanks to Lemma \ref{lemma:boundarybehaviour}, if
$U_2$ is not integrable at $0^+$ we can extend the map $C$ to a
continuous map defined on the whole space $(H, \|\cdot\|)$  defining ${C}\equiv
0$ on $\mathcal{D}(V)^c$. We set
$$G(\eta):= r\eta_0+ f\left(\eta\right)-{C}\left(\eta\right), \ \ \ \eta\in H,$$
and note that $G$ is continuous.

\smallskip
\emph{Local existence.} Let $\bar{\eta}\in{H}$ the initial datum for the equation.
We have to show the local existence of a solution of
\begin{equation*}
\begin{cases}
x'(t)=G\left((x(t),x(t+\xi)|_{\xi\in[-T,0]})\right),\\
x(0)=\bar{\eta}_0, \ x(s)=\bar{\eta}_1(s), \ s\in[-T,0),
\end{cases}
\end{equation*}
Since $G$ is continuous, there exists  $b>0$ such that
$m:=\sup_{\|\eta-\bar{\eta}\|^2\leq b} |G(\eta)|<+\infty .$ By
continuity of translations in $L^2(\mathbb{R};\mathbb{R})$ we can
find $a\in[0,T]$ such that
$$\int_{-T}^{-t}|\bar{\eta}_1(t+\xi)-\bar{\eta}_1(\xi)|^2d\xi\leq b/4, \ \ \ \ \forall t\in[0,a];$$
moreover, without loss of generality, we can suppose that
$\int_{-a}^0|\bar{\eta}_1(\xi)|^2d\xi\leq b/16$. Set
$$\displaystyle{\alpha:=\mbox{min}\left\{a,\,\frac{b}{2m},\, \frac{b}{16}\left(b+2|\bar{\eta}_0|^2\right)^{-1}\right\}}.
$$
Define
$$M:=\left\{x(\cdot)\in C\left([0,\alpha];\mathbb{R}\right) \ \big| \ |x(\cdot)-\bar{\eta}_0|^2\leq b/2\right\};$$
$M$ is a convex closed subset of the Banach space
$C([0,\alpha];\mathbb{R})$ endowed with the sup-norm. Define
$$x(t+\xi):=\bar{\eta}_1(t+\xi), \ \ \ \ \mbox{if} \ t+\xi\leq 0,$$
and observe that, for $t\in[0,\alpha]$, $x(\cdot)\in M$,
\begin{eqnarray*}
\int_{-t}^0|x(t+\xi)-\bar{\eta}_1(\xi)|^2d\xi&\leq& \int_{-t}^0\left(2|x(t+\xi)|^2+2|\bar{\eta}_1(\xi)|^2\right)d\xi\\
&\leq & 2\left[\int_{-t}^0\left(2\left(|x(t+\xi)-\bar{\eta}_0|\right)^2+2|\bar{\eta}_0|^2\right)d\xi+\int_{-t}^0 |\bar{\eta}_1(\xi)|^2d\xi\right]\\
&\leq
&2\left[2t\left(\frac{b}{2}+|\bar{\eta}_0|^2\right)+\frac{b}{16}\right]\leq
b/4
\end{eqnarray*}
So, for $t\in[0,\alpha]$, $x(\cdot)\in M$,  we have
\begin{eqnarray*}
\left\|\left(x(t),x(t+\xi)|_{\xi\in[-T,0]}\right)-\bar{\eta}\right\|^2&\leq&|x(t)-\bar{\eta}_0|^2
+\int_{-t}^0|x(t+\xi)-\bar{\eta}_1(\xi)|^2d\xi
+\int_{-T}^{-t}|\bar{\eta}_1(t+\xi)-\eta_1(\xi)|^2d\xi\\
&\leq & b/2+b/4+b/4=b.
\end{eqnarray*}
Define, for $t\in[0,\alpha]$, $x(\cdot)\in M$,
$$[\mathcal{J}x](t):=\bar{\eta}_0+ \int_0^t G\left(x(s),x(s+\xi)|_{\xi\in[-T,0]}\right)ds, \ \ \ t\in[0,\alpha].$$
We have
\begin{eqnarray*}
\Big|\, [\mathcal{J}x](t)-\eta_0\Big|&\leq &\int_0^t \left|\, G\left(x(s),x(s+\xi)|_{\xi\in[-T,0]}\right)\right|ds\\
&\leq & t m\leq b/2.
\end{eqnarray*}
Therefore we have proved that $\mathcal{J}$ maps the closed and
convex set $M$ in itself. We want to prove that $\mathcal{J}$ admits
a fixed point, i.e., by definition of $\mathcal{J}$, the solution we
are looking for.  By Schauder's Theorem it is enough to prove that
$\mathcal{J}$ is completely continuous, i.e. that
$\overline{\mathcal{J}(M)}$ is compact.
For any $x(\cdot)\in M$, we have  the estimate
$$\Big| \,[\mathcal{J}x](t)-[\mathcal{J}x](\bar{t})\Big|\leq\int_{t\wedge\bar{t}}^{t\vee\bar{t}}\left|\, G\left(x(s),x(s+\xi)|_{\xi\in[-T,0]}\right)\right|ds\leq m|t-\bar{t}|, \ \ \ \ t,\bar{t}\in[0,\alpha].$$
Therefore  $\mathcal{J}(M)$ is a uniformly bounded and
equicontinuous family in the space $C([0,\alpha];\mathbb{R})$. Thus,
by Ascoli-Arzel\a Theorem, $\overline{\mathcal{J}(M)}$ is compact.

\smallskip
\emph{Global existence.} Let $\eta\in H_{++}$ and let $x^*(\cdot)$ be  the solution of equation \eqref{rfe} defined on an interval $[0,\beta)$, $\beta>0$. Note that, by continuity of $f_0,\,{C}$, we have $x^*(\cdot)\in C^1([0,\beta);\mathbb{R})$.\\
Since ${C}(\cdot)\geq 0$, we have $x^*(\cdot)\leq x(\cdot;\eta,0)$; therefore $x^*(\cdot)$ is dominated from above on $[0,\beta)$ by $$\max_{t\in[0,\beta]} x(\cdot;\eta,0).$$\\
We want to show that it is also dominated from below in order to apply the extension argument. Let us suppose that $x^*(\bar{t})=0$ for some $\bar{t}\in[0,\beta)$. We want to show that this leads to a contradiction, so that, without loss of generality, we can suppose that
$$\bar{t}=\min\{t\in(0,\beta) \ | \ x^*(t)=0\}.$$
Therefore $x^*(\cdot)>0$ in a left neighborhood of $\bar{t}$. Since $f_0$ satisfies \eqref{f_0strict} and thanks to \eqref{a>0}-(ii), we must have $\frac{d}{dt}x^*(\bar{t})>0$, which contradicts $x^*(\cdot)>0$ in a left neighborhood of $\bar{t}$. Therefore we can say that $x^*(\cdot)>0$ on $[0,\beta)$, so that in particular $x^*(\cdot)$ is bounded from below by $0$ on $[0,\beta)$.
Therefore, arguing as in the classical extension theorems for ODE, we could show that we can extend $x^*(\cdot)$ to a solution defined on $[0,+\infty)$ and, again   by the same argument above, it will be $x^*(\cdot)>0$ on $[0,+\infty)$.\hfill$\square$

\section{Approximation results}\label{sec:appr}
In this section we obtain some approximation results which may be
used in order to produce $\varepsilon$-optimal controls for a wider
class of problems. Herein we assume that
\begin{equation}\label{continuasopra}
rx+f_0(x,0)\geq 0,  \ \ \ \forall x\geq 0,
\end{equation}
which implies in particular
\begin{equation}\label{conditionaltra}
\exists \delta>0 \ \mbox{such that} \ rx+f_0(x,0)\geq 0, \ \ \forall x\in (0,\delta].
\end{equation}
In fact, all the
results given below hold under \eqref{conditionaltra} as well. We
assume \eqref{continuasopra} only to simplify the proofs. Moreover
we incorporate the term $rx$ in the state equation within the term
$f_0$, so that consistently with \eqref{continuasopra} we assume
that
\begin{equation}\label{con:appr}
 f_0(x,0)\geq 0,  \ \ \ \forall x\geq 0.
\end{equation}

Recall that, for a control problem, an $\varepsilon$-optimal strategy is a strategy $\varepsilon$-near to optimality. Precisely, in our problem
\begin{definition}
 Let $\eta\in \mathcal{D}(V), \ \varepsilon >0$; an admissible control $c^\varepsilon(\cdot)\in\mathcal{C}(\eta)$ is said $\varepsilon$-optimal for the initial state $\eta$ if $J(\eta; c^\varepsilon(\cdot))>V(\eta)-\varepsilon$.     \hfill$\square$
\end{definition}

We will use the same concept for the control problems defined in the following.

\subsection{The case without utility on the state}
In the previous section we introduced an assumption of no
integrability of the utility function $U_2$. This was necessary in
order to ensure the existence of solutions for the closed loop
equation and the admissibility of the feedback strategy. This fact
is quite uncomfortable, because usually in consumption problems the
objective functional is given by a utility depending only on the
consumption variable, i.e. the case $U_2\equiv 0$ should be
considered. Of course we could take a $U_2$ heavily negative in a
right neighborhood of $0$ and equal to $0$ out of this neighborhood,
considering this as a forcing on the state constraint (states too
near to $0$ must be avoided). However we want to give here an
approximation procedure to partly treat also the case $U_2\equiv 0$,
giving a way to construct at least $\varepsilon$-optimal strategies
in this case.

So, let us consider a sequence of real functions $(U_2^n)$  such that
\begin{equation}\label{U_n}
U_2^n\uparrow 0, \ \ \ U_2^n \ \mbox{not integrable at} \ 0^+, \ \ \ U_2^n\equiv 0 \ \mbox{on} \ [1/n,+\infty).
\end{equation}
Let us denote by $J^n$ and $V^n$ respectively the objective
functionals and the value functions of the problems where the
utility on the state is given by $U_2^n$ and by $J^0$ and  $V^0$
respectively the objective functional and the value function of the
problem where the utility on the state disappears, i.e. $U_2\equiv
0$. It is immediate to see that monotonicity implies
\begin{equation}\label{convmon}
V^n\uparrow g\leq V^0.
\end{equation}
Thanks to the previous section, for any problem $V^n$,
$n\in\mathbb{N}$, we have an optimal feedback strategy
$c^*_n(\cdot)$.
\begin{lemma}\label{epsoptlemma}
Let $\eta\in\mathcal{D}(V^0)\subset H_+$. Then, for any $\varepsilon>0$, there exists an $\varepsilon$-optimal strategy $c^\varepsilon(\cdot)\in\mathcal{C}(\eta)$ for $V^0(\eta)$ such that
$$\inf_{t\in[0,+\infty)} x(t;\eta,c^{\varepsilon}(\cdot))>0.$$
\end{lemma}
\textbf{Proof.} Let $\varepsilon>0$ and take an $\varepsilon/2$-optimal control $c^{\varepsilon/2}(\cdot)\in\mathcal{C}(\eta)$ for $V^0(\eta)$. Let $M>T$ be such that
\begin{equation}\label{ghfd}
\frac{1}{\rho}e^{-\rho M} (\bar{U}_1-U_1(0))<\varepsilon/2.
\end{equation}
Define the control
$$c^\varepsilon(t):=\begin{cases}
c^{\varepsilon/2}(t), \ \ \ \ \ \ \ \ \ \mbox{for} \ t\in[0,M],\\
0, \ \ \ \ \ \ \ \ \ \ \ \ \ \ \mbox{for} \ t>M.
\end{cases}
$$
By Lemma \ref{comparison} we have
\begin{equation*}
x(\cdot;\eta,c^\varepsilon(\cdot))\geq x(\cdot;\eta,c^{\varepsilon/2}(\cdot))
\end{equation*}
and, by the assumption \eqref{con:appr} and since $c^\varepsilon(t)=0$ for $t\geq M$, it is not difficult to see that
$$x(t;\eta,c^\varepsilon(\cdot))\geq x(M;\eta,c^\varepsilon(\cdot)), \ \ \ \ \ \mbox{for} \ t\geq M,$$
so that
$$\inf_{t\in[0,+\infty)}x(t;\eta,c^\varepsilon(\cdot))=\inf_{t\in[0,M]} x(t;\eta,c^\varepsilon(\cdot))>0.$$
We claim that $c^\varepsilon(\cdot)$ is $\varepsilon$-optimal for $V^0(\eta)$, which yields the claim. Since $c^{\varepsilon/2}(\cdot)$ is $\varepsilon/2$-optimal for $V^0(\eta)$, taking also into account  \eqref{ghfd},  we get
\begin{multline*}
V^0(\eta)-\int_0^{+\infty} e^{-\rho t}U_1(c^\varepsilon(t))dt
=V^0(\eta)-\int_0^{+\infty}e^{-\rho t}U_1(c^{\varepsilon/2}(t))dt\\
+\int_M^{+\infty}e^{-\rho t}(U_1(c^{\varepsilon/2}(t))-U_1(0))dt<\varepsilon/2+\varepsilon/2=\varepsilon.
\end{multline*}
\hfill $\square$
\begin{proposition}\label{PR:nostateut}
Let $\eta \in \mathcal{D}(V^0)$ and $\varepsilon >0$. Then $V^n(\eta)\rightarrow V^0(\eta)$ and, when $n$ is large enough, $c^*_n(\cdot)$ is $\varepsilon$-optimal for $V^0(\eta)$.
\end{proposition}
\textbf{Proof.}
Let $\varepsilon>0$ and take an $\varepsilon$-optimal control $c^\varepsilon(\cdot)\in\mathcal{C}(\eta)$ for $V^0(\eta)$ such that (Lemma \ref{epsoptlemma})
$$m:=\inf_{t\in[0,+\infty)}x(t;\eta;c^\varepsilon(\cdot))>0.$$
Take $n\in\mathbb{N}$ such that $1/n<m$. Since $U_2^n\equiv 0$ on $[m,+\infty)$, we have
\begin{multline*}
V^0(\eta)-\varepsilon \leq  J(\eta;c^\varepsilon(\cdot))=
\int_0^{+\infty} e^{-\rho t} U_1(c^\varepsilon(t))dt\\= \int_0^{+\infty} e^{-\rho t} \left[ U_1(c^\varepsilon(t))+U_2^n(x(t;\eta,c^\varepsilon(\cdot)))\right]dt\\
=J^n(\eta,c^\varepsilon(\cdot)) \leq V^n(\eta)=J^n(\eta,c^*_n(\cdot))\leq J^0(\eta, c_n^*(\cdot)).$$
\end{multline*}
The latter inequality, toghether with \eqref{convmon}, proves both the claims.
\hfill $\square$

\subsection{The case with pointwise delay in the state
equation}\label{subs:apprBambi}

In this subsection we want to show that our problem is a good
approximation for growth models with time to build and concentrated
lag and discuss why our approach cannot work directly when the delay
is concentrated in a point. In this case the state equation is
\begin{equation}\label{eqstateak}
\begin{cases}
y'(t)=f_0\left(y(t),y\left(t-\frac{T}{2}\right)\right)-c(t),\\
y(0)=\eta_0, \ y(s)=\eta_1(s), \ s\in[-T,0).
\end{cases}
\end{equation}
It is possible to prove, as done in Theorem \ref{daprato}, that this
equation admits, for every $\eta\in H_+$, and for every $c(\cdot)\in
L^1_{loc} ([0,+\infty);\mathbb{R})$, a unique absolutely continuous
solution. We denote this  solution by $y(\cdot;\eta, c(\cdot))$. The
aim is  to maximize, over the set
\begin{equation}\label{defC^0}
\mathcal{C}^0_{ad}(\eta):=\{c(\cdot)\in L^1_{loc} ([0,+\infty);\mathbb{R}) \ | \ y(\cdot;\eta,c(\cdot))>0\},
\end{equation}
the functional
$$J_{0}(\eta, c(\cdot))=\int_0^{+\infty}e^{-\rho t}\left[U_1(c(\cdot))+U_2(y(t;\eta,c(\cdot)))\right]dt.$$
Denote by $V_0$ the associated value function. By monotonicity of $f_0$ we straightly get $H_{++}\subset\mathcal{D}(V_0)$.

Let us take a sequence $(a_k)_{k\in\mathbb{N}}\subset W^{1,2}_{-T}$  such that
\begin{equation}\label{defa_k}
a_k(-T)=0, \ \  \ \|a_k\|_{L^2_{-T}}\leq 1, \ \ \mbox{\eqref{a>0}-(ii) holds true} \ \forall a_k, \ \ a_k\stackrel{*}{\rightharpoonup}\delta_{-T/2} \ \mbox{in} \ \left(C([-T,0];\mathbb{R})\right)^*,
\end{equation}
where $\delta_{-T/2}$ is the Dirac measure concentrated at $-T/2$. We denote by $x_k(\cdot;\eta,c(\cdot))$ the unique solution of \eqref{eqstate} where $a(\cdot)$ is replaced by $a_k(\cdot)$.
\begin{proposition}\label{comparak}
Let $\eta\in H_+$, $c(\cdot)\in L^1_{loc}([0,+\infty);\mathbb{R})$ and set  $y(\cdot):=y(\cdot;\eta,c(\cdot))$,  $x_k(\cdot):=x_k(\cdot;\eta,c(\cdot))$. Then there exists a continuous and increasing function $h$ such that $h(0)=0$ and
\begin{equation}\label{estimatestate}
\sup_{s\in[0,t]}|x_k(s)-y(s)|\leq h(t) u_k(t), \ \ \ \ t\in[0,+\infty),
\end{equation}
where $u_k(t)\rightarrow 0$, as $k\rightarrow\infty$, unifromly on
bounded sets.
\end{proposition}
\textbf{Proof.}
Note that
\begin{equation}\label{stimasua_k}
\|a_k\|_{\left(C([-T,0];\mathbb{R})\right)^*}=\sup_{\|f\|_{\infty}=1}\left|\int_{-T}^0 a_k(\xi)f(\xi)d\xi\right|\leq \int_{-T}^0|a_k(\xi)|d\xi\leq \|a_k\|_{L^2_{-T}}\cdot T^{1/2}\leq T^{1/2}.
\end{equation}
Let $t\geq 0$; we have, for any $\zeta\in[0,t]$,
\begin{eqnarray}\label{opiuy}
|x_k(\zeta)-y(\zeta)|&=&\int_0^\zeta \left[f_0\left(x_k(s),\int_{-T}^0a_k(\xi)x_k(s+\xi)d\xi\right)-f_0\left(y(s),y\left(s-\frac{T}{2}\right)\right)\right]ds\nonumber\\
&\leq &C_{f_0} \Bigg[\int_0^t |x_k(s)-y(s)|ds+\int_0^t \left|\int_{-T}^0a_k(\xi)(x_k(s+\xi)-y(s+\xi))d\xi\right| ds\nonumber\\
&&+\int_0^t \left|\int_{-T}^0a_k(\xi)y(s+\xi)d\xi- y\left(s-\frac{T}{2}\right)\right|ds\Bigg].
\end{eqnarray}
Call
$$g_k(t):= \sup_{s\in[-T,t]} |x_k(s)-y(s)|=\sup_{s\in[0,t]} |x_k(s)-y(s)|,$$
and set, for $s\in [0,t]$,
$$u_k(s):= C_{f_0}\int_0^s \left|\int_{-T}^0a_k(\xi)y(r+\xi)d\xi- y\left(r-\frac{T}{2}\right)\right|dr.$$
Note that, for every $s\in[0,t]$, the function $[-T,0]\ni \xi\mapsto x_k(s+\xi)-y(s+\xi)$ is  continuous, therefore thanks to \eqref{stimasua_k} we can write from \eqref{opiuy}
$$g_k(t)\leq C_{f_0}\left[ \int_0^t g_k(s)ds +T^{1/2} \int_0^t g_k(s)ds+ u_k(t)\right].$$
Therefore, setting $K:=C_{f_0}(1+T^{1/2})$, we get by Gronwall's Lemma
\begin{equation}\label{plkj}
g_k(t)\leq u_k(t)+Kte^{Kt}u_k(t)=:h(t)u_k(t).
\end{equation}
Note that, since $a_k\stackrel{*}{\rightharpoonup}\delta_{-T/2}$ in $\left(C([-T,0];\mathbb{R})\right)^*$, we have the pointwise convergence
$$\int_{-T}^0a_k(\xi)y(s+\xi)d\xi\longrightarrow y\left(s-\frac{T}{2}\right), \ \ \ s\in [0,t];$$
moreover
$$
\left|\int_{-T}^0 a_k(\xi)y(s+\xi)d\xi\right|\leq \|a_k\|_{L^2_{-T}}\cdot \left\|\,y(s+\xi)|_{\xi\in[-T,0]}\right\|_{L^2_{-T}}\leq C_{\eta,c(\cdot)}<+\infty, \ \ \ \ \forall s\in[0,t],$$
where the last inequality follows from the fact that the function $[0,t]\rightarrow L^2_{-T}$, $s\mapsto y(s+\xi)|_{\xi\in[-T,0]}$ is continuous.
Therefore we have by dominated convergence $u_k(t)\rightarrow 0$. By \eqref{plkj} we get \eqref{estimatestate}.\hfill$\square$\\

\noindent
For  $k\in\mathbb{N}$, $\eta\in H_{+},$ let
\begin{equation}\label{defC^k}
\mathcal{C}^k_{ad}(\eta):=\{c(\cdot)\in L^1_{loc}([0,+\infty);\mathbb{R}) \ | \ x_k(\cdot;\eta,c(\cdot))>0\},
\end{equation}
Consider the problem of maximizing over $\mathcal{C}^k_{ad}(\eta)$ the functional
$$J_{k}(\eta, c(\cdot)):=\int_0^{+\infty} e^{-\rho t} \left[U_1(c(t))+U_2(x_k(t;\eta,c(\cdot)))\right]dt$$
and denote by $V_k$ the associated value function.
Note that, since we have assumed \eqref{con:appr}, straightly we get  $H_{++}\subset \mathcal{D}(V_k)$ for every $k\in\mathbb{N}$.
Thanks to the previous section we have a sequence of optimal feedback strategies for the sequence of problems $(V_k(\eta))_{k\in\mathbb{N}}$, in the sense that we have a sequence $(c^*_k(\cdot))_{k\in\mathbb{N}}$ of feedback controls such that $c_k^*(\cdot)\in\mathcal{C}^k_{ad}(\eta)$ for every $k\in\mathbb{N}$ and
$$J_k(\eta;c^*_k(\cdot))=\sup_{c(\cdot)\in\mathcal{C}^k_{ad}(\eta)}J_k(\eta; c(\cdot))=:V_k(\eta), \ \ \ \ \ \ \  \forall k\in\mathbb{N}.$$
\begin{lemma}\label{lemmaepsok}
Let $\eta\in H_{++}$.
\begin{itemize}
\item For any $\varepsilon>0$ there exists an $\varepsilon$-optimal
strategy $c^{\varepsilon}(\cdot)\in\mathcal{C}_{ad}^0(\eta)$ for the problem $V_{0}(\eta)$ such that
$$\inf_{t\in[0,+\infty)} y(t;\eta,c^\varepsilon(\cdot))>0.$$
\item Assume that
\begin{equation}\label{limxux}
\lim_{x\rightarrow 0^+}\left[x\,U_2(x)\right]=-\infty.
\end{equation}
Then, {for any $\varepsilon>0$ there exists $\nu>0$ such that for any $k\in\mathbb{N}$ there exists} an  $\varepsilon$-optimal control $c^\varepsilon_k(\cdot)\in\mathcal{C}^k_{ad}(\eta)$ for the problem $V_{k}(\eta)$ such that
\begin{equation}\label{bebe}
\inf_{t\in[0,+\infty)} x_k(t;\eta,c^\varepsilon_{k}(\cdot))\geq \nu.
\end{equation}
\end{itemize}
\end{lemma}
\textbf{Proof.} (i) Let $\varepsilon>0$ and take an $\varepsilon/2$-optimal control $c^{\varepsilon/2}(\cdot)\in\mathcal{C}^0_{ad}(\eta)$ for the problem $V_{0}(\eta)$. Take $M>0$ large enough to satisfy
\begin{equation}\label{assM}
\frac{1}{\rho} e^{-\rho M} (\bar{U}_1- U_1(0))<\varepsilon/2.
\end{equation}
Define the control
$$c^\varepsilon(t):=\begin{cases}c^{\varepsilon/2} (t), \ \ \ \ \ \ \ \ \mbox{for} \ t\in [0,M],\\
0, \ \ \ \ \ \ \ \ \ \ \ \ \ \ \ \mbox{for} \ t>M.
\end{cases}
$$
A comparison criterion like the one proved in Lemma \ref{comparison} can be proved also for  equation \eqref{eqstateak}. Therefore we have
\begin{equation}\label{fdsae}
y(\cdot;\eta,c^\varepsilon(\cdot))\geq y(\cdot;\eta, c^{\varepsilon/2}(\cdot))
\end{equation}
and, since we have assumed \eqref{con:appr},
$$\inf_{t\in[0,+\infty)} y(t;\eta,c^\varepsilon(\cdot))=\inf_{t\in[0,M]} y(t;\eta,c^\varepsilon(\cdot)).$$
We claim that $y(\cdot;\eta,c^\varepsilon(\cdot))$ is $\varepsilon$-optimal for $V_{0}$, which yields the claim. Since $c^{\varepsilon/2}(\cdot)$ is $\varepsilon/2$-optimal for $V_{0}(\eta)$, taking also into account \eqref{assM} and \eqref{fdsae},
\begin{multline*}
V_{0}(\eta)-\int_0^{+\infty} e^{-\rho t} \left[U_1(c^\varepsilon(t))+U_2(y(t;\eta, c^\varepsilon(\cdot)))\right]dt\\
=V_{0} (\eta)-\int_0^{+\infty} e^{-\rho t}(U_1(c^{\varepsilon/2}(t))+U_2(y(t;\eta,c^{\varepsilon/2}(\cdot)))dt\\+\int_M^{+\infty} e^{-\rho t} (U_1(c^{\varepsilon/2}(t))-U_1(0))dt+\int_M^{+\infty}e^{-\rho t} \left[U_2(y(t;\eta,c^{\varepsilon/2}(\cdot))-U_2(y(t;\eta,c^\varepsilon (t)))\right]dt\\
<\varepsilon/2+\varepsilon/2=\varepsilon.
\end{multline*}

(ii) Due to \eqref{con:appr} we have $x_k(\cdot;\eta,0)\geq \eta_0$ for every $k\in \mathbb{N}$. Let
$$j_0:=\frac{ U_1(0)+U_2(\eta_0)}{\rho}.$$
Then we have $V_{k}(\eta)\geq J_{k}(\eta, 0)\geq  j_0$ for every $k\in\mathbb{N}$. Take $M>0$ large enough to satisfy
\begin{equation}\label{assM2}
\frac{1}{\rho} e^{-\rho M} (\bar{U}_1- U_1(0))<\varepsilon.
\end{equation}
Arguing as done to get \eqref{ss} and taking into account the comparison criterion, we can find $C_M>0$ such that, for every $k\in\mathbb{N}$ and for every $c(\cdot)\in\mathcal{C}^k_{ad}(\eta)$, we have for all $s\in[0,M]$ and for all $t\in[s,M+1]$,
\begin{equation}\label{rte}
x_k(t;\eta,c(\cdot))\leq x_k(s;\eta, c(\cdot))(1+C_M(t-s))+C_M(t-s).
\end{equation}
Now take $\nu>0$ small enough to have
\begin{equation}\label{nusmall}
\mbox{(i)} \ \ \nu<1, \ \ \ \mbox{(ii)} \ \ \frac{\nu}{2C_M}<1, \ \ \ \mbox{(iii)} \ \ \frac{\nu}{2C_M}U_2(2\nu) e^{-\rho (M+1)}<j_0-\frac{\bar{U}_1+\bar{U}_2}{\rho}-1<0.
\end{equation}
For $k\in\mathbb{N}$, thanks to the previous section we have optimal strategies in feedback form  $c^{*}_{k}(\cdot)\in\mathcal{C}^k_{ad}(\eta)$ for $V_{k}$; we claim that $x_k(t;\eta,c^{*}_{k}(\cdot))>\nu$ for $t\in[0,M]$ for every $k\in\mathbb{N}$. Indeed suppose by contradiction that for some $t_0\in [0,M]$ we have $x_k(t_0;\eta,c^{*}_{k}(\cdot))=\nu$; then by \eqref{rte} and \eqref{nusmall}-(i),(ii) we get that
$$x_k(t;\eta, c^{*}_{k}(\cdot))\leq 2\nu, \ \ \ \ \ \mbox{for} \ t\in\left[t_0,t_0+\frac{\nu}{2C_M}\right].$$
Therefore, by \eqref{nusmall}-(iii),
$$\int_{t_0}^{t_0+\frac{\nu}{2C_M}}e^{-\rho t}  U_2^n(x_k(t;\eta,c^{*}_{k}(\cdot)))dt\leq j_0-\frac{\bar{U}_1+\bar{U}_2}{\rho} -1.$$
This shows that
$$J_{k}(\eta, c^{*}_{k}(\cdot))\leq j_0-1\leq V_{k}(\eta)-1.$$
This fact contradicts the optimality of $c^{*}_{k}(\cdot)$. Therefore we have proved that for the choice of $\nu$ given by \eqref{nusmall} we have
$$x_k(t;\eta, c^{*}_{k}(\cdot))>\nu, \ \ \ \ \ \mbox{for} \ t\in[0,M].$$
We can continue the strategy $c^{*}_{k}(\cdot)$ after $M$ taking the null strategy, i.e. defining the strategy
\begin{equation}\label{rottopalle}
c^\varepsilon_{k}(\cdot):=\begin{cases}
c^{*}_{k}(t), \ \ \ \ \ \ \ \mbox{for} \ t\in [0,M],\\
0, \ \ \ \  \ \ \ \ \ \ \ \ \mbox{for} \ t>M.
\end{cases}
\end{equation}
Then by \eqref{con:appr} we have $x_k(\cdot;\eta, c^\varepsilon_{k}(\cdot))>\nu$ for every $k\in\mathbb{N}$. We claim that $c^\varepsilon_{k}(\cdot)$ is $\varepsilon$-optimal for $V_{k}(\eta)$ for every $k\in\mathbb{N}$, which proves the claim.
Indeed, taking into account the comparison criterion and \eqref{assM2} for the inequality in the following,
\begin{multline*}
V_{k}(\eta)-\int_0^{+\infty} e^{-\rho t} (U_1(c^\varepsilon_{k}(t))+U_2(x_k(t;\eta, c^\varepsilon_{k}(\cdot)))dt\\
=V_{k} (\eta)-\int_0^{+\infty} e^{-\rho t}(U_1(c^{*}_{k}(t))+U_2(x_k(t;\eta,c^{*}_{k}(\cdot)))dt\\+\int_M^{+\infty} e^{-\rho t} (U_1(c^{*}_{k}(t))-U_1(0))dt+\int_M^{+\infty}e^{-\rho t} (U_2(x_k(t;\eta,c^{*}_{k}(\cdot))-U_2(x_k(t;\eta,c^\varepsilon_{k} (t))))dt<\varepsilon.
\end{multline*}
\hfill$\square$

\begin{proposition}\label{propV_n,0}
Let $\eta\in H_{++}$ and suppose that \eqref{limxux} holds true. We
have $V_k(\eta)\rightarrow V_0(\eta)$, as $k\rightarrow\infty$.
Moreover for every $\varepsilon>0$  we can find a constant
$M_\varepsilon$ and a $k_\varepsilon$ such that the strategy
($c_k^*$ is the optimal feedback strategy for the problem of $V_k$)
\begin{equation}\label{mnbvnew}
c_{k_\varepsilon,{M}_\varepsilon}(t):=\begin{cases}
c^*_{k_\varepsilon}(t), \ \ \ \ \ \ \mbox{for} \ t\in[0,{M}_\varepsilon],\\
0, \ \ \ \ \ \ \ \ \ \ \ \mbox{for} \ t>{M}_\varepsilon.
\end{cases}
\end{equation}
is $\varepsilon$-optimal strategy for the problem $V_0(\eta)$.
\end{proposition}
\textbf{Proof.}
(i) Here we show that
\begin{equation}\label{liminfV_n}
\liminf_{k\rightarrow \infty} V_k(\eta)\geq V_0(\eta),
\end{equation}
Let $\varepsilon>0$ and let
$c^\varepsilon(\cdot)\in\mathcal{C}_{ad}^0(\eta)$ be an
$\varepsilon$-optimal strategy for the problem $V_{0}(\eta)$. Thanks
to Lemma \ref{lemmaepsok}-(i) we can suppose without loss of
generality $2\nu_1:=\inf_{t\in[0,+\infty)}y(t;\eta,
c^\varepsilon(\cdot))>0$. The function $U_1$ is uniformly continuous
on $[0,+\infty)$ and the function $U_2$ is uniformly continuous on
$[\nu_1,+\infty)$. Let $\omega_{\nu_1}$ be a modulus of uniform
continuity for both these functions. Take $M>0$ such that
\begin{equation}\label{hjk}
- \frac{1}{\rho}e^{-\rho M}(\bar{U}_1+\bar{U}_2)-\frac{1-e^{-\rho M}}{\rho} \omega_{\nu_1}(\nu_1)+\frac{1}{\rho} e^{-\rho M} (U_1(0)+U_2(\nu_1))
\geq -\varepsilon,
\end{equation}
Define
$$c_M^\varepsilon(t):=\begin{cases}
c^\varepsilon(t), \ \ \ \ \ \ \mbox{for} \ t\in[0,M],\\
0,\ \  \ \ \ \ \ \ \ \ \ \mbox{for} \ t>M.
\end{cases}
$$
Let $k_{M}$ be such that
\begin{equation}\label{poop}
h(M)u_k(M)< \nu_1, \ \ \ \forall k\geq k_{M},
\end{equation}
where $u_k$ and $h$ are  the functions appearing in
\eqref{estimatestate}. Then, thanks to Proposition \ref{comparak}
and to the  monotonicity property of $f_0$, it is straightforward to
see that $x_k(t;\eta,c_M^\varepsilon(\cdot))\geq \nu_1>0$, so that
in particular $c^\varepsilon_M(\cdot)\in \mathcal{C}_{ad}^k(\eta)$
for all $k\geq k_{M}$. For all $k\geq k_{M}$, we have, thanks to
Proposition \ref{comparak} and by definition of $k_{M}$,
\begin{multline*}
\int_0^{+\infty} e^{-\rho t} [U_1(c_M^\varepsilon(t))+U_2(x_k(t;\eta,c_M^\varepsilon(\cdot)))]dt\\=\int_0^{M} e^{-\rho t} [U_1(c_M^\varepsilon(t))+U_2(x_k(t;\eta,c_M^\varepsilon(\cdot)))]dt+\int_M^{+\infty} e^{-\rho t} [U_1(c_M^\varepsilon(t))+U_2(x_k(t;\eta,c_M^\varepsilon(\cdot)))]dt\\
\geq  \int_0^M   e^{-\rho t} [U_1(c^\varepsilon(t))+U_2(y(t;\eta, c^\varepsilon(\cdot)))]dt-\frac{1-e^{-\rho M}}{\rho} \omega_{\nu_1}(\nu_1)
+\frac{1}{\rho} e^{-\rho M} (U_1(0)+U_2(\nu_1))\\
\geq  \int_0^{+\infty}   e^{-\rho t} [U_1(c^\varepsilon(t))+U_2(y(t;\eta, c^\varepsilon(\cdot)))]dt\\- \frac{1}{\rho}e^{-\rho M}(\bar{U}_1+\bar{U}_2)-\frac{1-e^{-\rho M}}{\rho} \omega_{\nu_1}(\nu_1)+\frac{1}{\rho} e^{-\rho M} (U_1(0)+U_2(\nu_1)).
\end{multline*}
so that by \eqref{hjk}
\begin{equation}\label{fhg}
V_{k}(\eta)\geq V_{0}(\eta)-2\varepsilon,
\end{equation}
which shows \eqref{liminfV_n}.

\smallskip
(ii) Now we show that
\begin{equation}\label{wlimsup}
\limsup_{k\rightarrow \infty} V_k(\eta)\leq V_0(\eta).
\end{equation}
Let $\varepsilon>0$; thanks to Lemma \ref{lemmaepsok}-(ii) we can construct a sequence $(c^\varepsilon_k(\cdot))_{k\in\mathbb{N}}$, $c^\varepsilon_k(\cdot)\in\mathcal{C}^k_{ad}(\eta)$ for every $k\in\mathbb{N}$,  of $\varepsilon$-optimal controls for the sequence of problems $(V_k(\eta))_{k\in\mathbb{N}}$ such that
$$2\nu_2:=\inf_{k\in\mathbb{N}} \inf_{t\in[0,+\infty)} x_k(t;\eta,c^\varepsilon_k(\cdot))>0.$$
Let $\omega_{\nu_2}$ be a modulus of uniform continuity for $U_1$ on $[0,+\infty)$ and for $U_2$ on $[\nu_2,+\infty)$.
Take $\tilde{M}>0$ such that
\begin{equation}\label{tildeM}
- \frac{1}{\rho}e^{-\rho \tilde{M}}(\bar{U}_1+\bar{U}_2)-\frac{1-e^{-\rho \tilde{M}}}{\rho} \omega_{\nu_2}(\nu_2)+\frac{1}{\rho} e^{-\rho \tilde{M}} (U_1(0)+U_2(\nu_2))>-\varepsilon
\end{equation}
and define the controls
\begin{equation}\label{mnbv}
c_{k,\tilde{M}}^\varepsilon(t):=\begin{cases}
c^\varepsilon_k(t), \ \ \ \ \ \ \mbox{for} \ t\in[0,\tilde{M}],\\
0, \ \ \ \ \ \ \ \ \ \ \ \mbox{for} \ t>\tilde{M}.
\end{cases}
\end{equation}
As before we can find $k_{\tilde{M}}$ such that we have
\begin{equation}\label{poop2}
h(\tilde{M})u_k(\tilde{M})< \nu_2, \ \ \ \forall k\geq k_{\tilde{M}},
\end{equation}
In this case we have $y(\cdot;\eta,c^\varepsilon_k(\cdot))\geq \nu_2$ and
\begin{multline}\label{pkk}
\int_0^{+\infty} e^{-\rho t} [U_1(c_{k,\tilde{M}}^\varepsilon(t))+U_2(y(t;\eta,c_{k,\tilde{M}}^\varepsilon(\cdot)))]dt\\=
\int_0^{\tilde{M}} e^{-\rho t} [U_1(c_{k,\tilde{M}}^\varepsilon(t))+U_2(y(t;\eta,c_{k,\tilde{M}}^\varepsilon(\cdot)))]dt+\int_{\tilde{M}}^{+\infty} e^{-\rho t} [U_1(c_{k,\tilde{M}}^\varepsilon(t))+U_2(y(t;\eta,c_{k,\tilde{M}}^\varepsilon(\cdot)))]dt\\
\geq  \int_0^{\tilde{M}}   e^{-\rho t} [U_1(c_k^\varepsilon(t))+U_2(x_k(t;\eta, c^\varepsilon_k(\cdot)))]dt-\frac{1-e^{-\rho \tilde{M}}}{\rho} \omega_{\nu_2}(\nu_2)+\frac{1}{\rho} e^{-\rho \tilde{M}} (U_1(0)+U_2(\nu_2))\\
\geq \int_0^{+\infty}   e^{-\rho t} [U_1(c_k^\varepsilon(t))+U_2(x_k(t;\eta, c^\varepsilon_k(\cdot)))]dt\\- \frac{1}{\rho}e^{-\rho \tilde{M}}(\bar{U}_1+\bar{U}_2)-\frac{1-e^{-\rho \tilde{M}}}{\rho} \omega_{\nu_2}(\nu_2)+\frac{1}{\rho} e^{-\rho \tilde{M}} (U_1(0)+U_2(\nu_2)).
\end{multline}
By \eqref{tildeM},
we get, for $k\geq k_{\tilde{M}}$,
$$V_0(\eta)\geq V_k(\eta)-2\varepsilon,$$
which proves \eqref{wlimsup}.

\smallskip
(iii) The procedure of construction of $c^\varepsilon_{k,M}$ in (ii)
yields $\varepsilon$-optimal controls for the limit problem
$V_0(\eta)$. Indeed, starting from $\varepsilon>0$, we can compute
$\nu_1,\nu_2,M,\tilde{M}$ {depending on $\varepsilon$} such that \eqref{hjk}  and \eqref{tildeM}
hold true.  Then, if $(a_k)_{k\in\mathbb{N}}$ is chosen in a clever
way, for example if $(a_k)_{k\in\mathbb{N}}$ is a sequence of
gaussian densities, we can compute $k_{M}, k_{\tilde{M}}$  such that
\eqref{poop}-\eqref{poop2} hold true. Thanks to \eqref{pkk} and
\eqref{fhg}, for every $k\geq k_{M}\vee k_{\tilde{M}}$ the controls
$c^\varepsilon_{k,\tilde{M}}(\cdot)$ defined in \eqref{mnbv} are
$4\varepsilon$-optimal for the limit problem $V^0(\eta)$. Replacing
$\varepsilon$ with $\varepsilon/4$ we get {the controls in \eqref{mnbvnew}}. \hfill
$\square$

\begin{remark}\label{rm:noBambi}
When the delay is concentrated in a point in a linear way, we could
tempted to insert the delay term in the infinitesimal generator $A$
and try to work as done in Section \ref{sec:infinite}. Unfortunately
this is not possible. Indeed consider this simple case:
\begin{equation*}
\begin{cases}
y'(t)=ry(t)+y\left(t-T\right),\\
y(0)=\eta_0, \ y(s)=\eta_1(s), \ s\in[-T,0),
\end{cases}
\end{equation*}
In this case we can define
$$
A:\mathcal{D}(A)\subset H\longrightarrow H, \qquad
(\eta_0,\eta_1(\cdot))\longmapsto  (r
\eta_0+\eta_1(-T),\eta_1'(\cdot)).
$$
where again
$$\mathcal{D}(A):=\{ \eta\in H \ | \ \eta_1(\cdot)\in W^{1,2}([-T,0]; \mathbb{R}), \ \eta_1(0)=\eta_0\}.$$
The inverse of $A$ is the operator
$$
A^{-1}:(H,\|\cdot\|)\longrightarrow (\mathcal{D}(A),\|\cdot\|)
\qquad(\eta_0,\eta_1(\cdot))\longmapsto
\displaystyle{\left(\frac{\eta_0-c}{r}, \ c+\int_{-T}^\cdot
\eta_1(\xi)d\xi\right),}
$$
where
$$c=\frac{1}{r+1}\,\eta_0-\frac{r}{r+1}\,\int_{-T}^0\eta_1(\xi)d\xi.$$
In this case we would have the first part of Lemma \ref{estimateA},
but not the second part, because it is not possible to control
$|\eta_0|$ by $\|\eta\|_{-1}$. Indeed take for example $r$ such that
$\frac{1-r}{1+r}=\frac{1}{2}$, and $(\eta^n)_{n\in\mathbb{N}}\subset
H$ such that
$$\eta_0^n=1/2, \ \ \ \int_{-T}^0\eta_1^n(\xi)d\xi=1, \ \ \ n\in\mathbb{N}.$$
We would have $c=1/2$, so that $\left|\frac{\eta_0^n-c}{r}\right|=0$. Moreover we can choose $\eta_1^n$ such that, when $n\rightarrow \infty$,
$$\int_{-T}^0\left|\,\frac{1}{2}+\int_{-T}^s\eta_1^n(\xi)d\xi\right|^2ds\longrightarrow 0.$$
Therefore we would have $|\eta_0^n|=1/2$ and $\|\eta^n\|_{-1}\rightarrow 0$. This shows that the second part of Lemma \ref{estimateA} does not hold. Once this part does not hold, then everything in the following argument breaks down.\hfill$\blacksquare$
\end{remark}

\subsection{The case with pointwise delay in the state equation and without
utility on the state}\label{subs:apprBambibis}

Now we want to approximate the problem of optimizing, for $\eta\in
H_{++}$,
$$J_{0}^{0}(\eta,c(\cdot)):=\int_{0}^{+\infty} e^{-\rho t}
U_1(c(t))dt,$$ over the set \eqref{defC^0}, where
$y(\cdot;\eta,c(\cdot))$ follows the dynamics given by
\eqref{eqstateak}. Let us denote by $V_0^0$ the corresponding value
function and let us take a sequence of real functions $(U_2^n)$ as
in \eqref{U_n}, but with the assumption of no integrability at $0^+$
replaced by the stronger assumption
$$\lim_{x\rightarrow 0^+} x\,U_2^n(x)=-\infty, \ \ \ \ \ \ \forall n\in\mathbb{N}.$$
Fix  $n\in\mathbb{N}$ and consider the sequence of functions
$(a_k)_{k\in\mathbb{N}}$ defined in \eqref{defa_k}. For
$k\in\mathbb{N}$, consider  the problem of maximizing over the set
$\mathcal{C}^k_{ad}(\eta)$ defined in \eqref{defC^k} the functional
$$J_k^n(\eta, c(\cdot)):=\int_0^{+\infty} e^{-\rho t} (U_1(c(t))+U_2^n(x_k(t;\eta,c(\cdot)))dt,$$
where $x_k(\cdot;\eta, c(\cdot))$ follows the dynamics given by \eqref{eqstate} when $a(\cdot)$ is replaced by $a_k(\cdot)$,
and denote by $V_k^n$ the associated value function.

Moreover, for $k\in\mathbb{N}$,  consider the problem of maximizing
$$J_{k}^{0}(\eta,c(\cdot)):=\int_{0}^{+\infty} e^{-\rho t} U_1(c(t))dt,$$
over $\mathcal{C}^k_{ad}(\eta)$ and denote by $V^0_k$ the associated value function.

Finally consider the problem of maximizing over the set
 $\mathcal{C}^0_{ad}(\eta)$ the functional
$$J_0^n(\eta, c(\cdot)):=\int_0^{+\infty} e^{-\rho t} (U_1(c(t))+U_2^n(y(t;\eta,c(\cdot)))dt,$$
and denote by $V_{0}^n$ the associated value function.

\smallskip
For fixed $n\in\mathbb{N}$, the problems $V_k^n$ approximate, when $k\rightarrow \infty$, the problem $V_0^n$ in the sense of Proposition \ref{propV_n,0}, i.e. we are able to produce
$k_{\varepsilon,n}, M_{\varepsilon,n}$  large enough  to make the strategy
$c_{k_{\varepsilon,n}, {M}_{\varepsilon,n}}(\cdot)$ defined as in \eqref{mnbvnew} admissible and $\varepsilon$-optimal for the problem $V_0^n(\eta)$.
\begin{proposition}\label{PR:both}
Let $\eta\in H_{++}$, let $k_{\varepsilon,n}, M_{\varepsilon,n}$,
$c_{k_{\varepsilon,n}, {M}_{\varepsilon,n}}(\cdot)$ as above. For every $\varepsilon>0$ we can find $n_\varepsilon$ such that
\begin{equation}\label{ultimoclaim}
\lim_{\varepsilon\downarrow 0}V_{k_{\varepsilon,n_\varepsilon}}^{n_\varepsilon}(\eta)= V_0^0(\eta).
\end{equation}
Moreover the controls $c_{k_{\varepsilon,n_\varepsilon}, {M}_{\varepsilon,n_\varepsilon}}(\cdot)$ defined as in \eqref{mnbvnew} are admissible and $3\varepsilon$-optimal for the problem $V_0^0(\eta)$.
\end{proposition}
\textbf{Proof.}
Let $\varepsilon>0$ and consider the strategies  $c^\varepsilon(\cdot)$ and   $c^\varepsilon_{M}(\cdot)$ defined as in the part (i) of the proof of Proposition \ref{propV_n,0}. Notice that actually $M=M(\varepsilon,n)=:M_n^\varepsilon$.
Notice also that by definition of $k_{\varepsilon,n}$, $M_n^\varepsilon$ we have $x_{k_{\varepsilon,n}}(\cdot;\eta, c_{M^\varepsilon_n}^\varepsilon(\cdot))\geq \nu_1$ and that \eqref{assM} in particular implies
$$\frac{1}{\rho}e^{-\rho M_{n}^\varepsilon} (\bar{U}_1-U_1(0))\leq \varepsilon.$$
Take $n_\varepsilon\in\mathbb{N}$ such that $1/n_\varepsilon<\nu_1$ (notice that $\nu_1$ depends on $\varepsilon$ and does not depend on $n$).  Then, since $U_2^{n_\varepsilon}\equiv 0$ on $[\nu_1,+\infty)$, we can write
\begin{multline}\label{ptet}
V_0^0(\eta)-\varepsilon \leq  J^0_0(\eta;c^\varepsilon(\cdot))=
\int_0^{+\infty} e^{-\rho t} U_1(c^\varepsilon(t))dt\\= \int_0^{+\infty} e^{-\rho t} \left[ U_1(c^\varepsilon(t))+U_2^{n_\varepsilon}(x_{k_{\varepsilon,n_\varepsilon}}(t;\eta,c_{M_{n_\varepsilon}^\varepsilon}^\varepsilon(\cdot)))\right]dt\\
\leq J^{n_\varepsilon}_{k_{\varepsilon,n_\varepsilon}}(\eta,c_{M^\varepsilon_{n_\varepsilon}}^\varepsilon(\cdot)) +\frac{1}{\rho}e^{-\rho M^\varepsilon_{n_\varepsilon}} (\bar{U}_1-U_1(0)) \leq V^{n_\varepsilon}_{k_{\varepsilon,n_\varepsilon}}(\eta)+\varepsilon,$$
\end{multline}
so that
\begin{equation}\label{liminf:ultimo}
\liminf_{\varepsilon\downarrow 0}V^{n_\varepsilon}_{k_{\varepsilon,n_\varepsilon}}(\eta)\geq V^0_0(\eta).
\end{equation}

\smallskip
On the other hand the strategies  $c_{k_{\varepsilon,n},M_{\varepsilon,n}}(\cdot)$ defined in \eqref{mnbvnew} are admissible   for the problem $V^0_0$ (since the state equation related to $V_0^n$ and to $V^0_0$ is the same) and $c_{k_{\varepsilon,n},{M}_{\varepsilon,n}}(\cdot)$ is $\varepsilon$-optimal for $V_{k_{\varepsilon,n}}^n(\eta)$ for every $n\in\mathbb{N}$. Therefore
 \begin{equation}\label{puyt}
V_{k_{\varepsilon,n_\varepsilon}}^{n_\varepsilon}(\eta)-\varepsilon \leq  J_{k_{\varepsilon,n_\varepsilon}}^{n_\varepsilon}(\eta;c_{k_{\varepsilon,n_\varepsilon},{M}_{\varepsilon,n_\varepsilon}}(\cdot))\leq J_{k_{\varepsilon,n_\varepsilon}}^0(\eta;c_{k_{\varepsilon,n_\varepsilon},{M}_{\varepsilon,n_\varepsilon}}(\cdot)) = J_0^0(\eta;c_{k_{\varepsilon,n_\varepsilon},{M}_{\varepsilon,n_\varepsilon}}(\cdot))\leq V^0_0(\eta),
\end{equation}
which shows
\begin{equation}\label{limsup:ultimo}
\limsup_{\varepsilon\downarrow 0}V^n_{k_{\varepsilon,n}}(\eta)\leq V^0_0(\eta).
\end{equation}

\smallskip
Combining \eqref{liminf:ultimo} and \eqref{limsup:ultimo} we get \eqref{ultimoclaim}. Combining \eqref{ptet} and \eqref{puyt} we get
$$V^0_0(\eta)\leq J_0^0(\eta;c_{k_{\varepsilon,n_\varepsilon},{M}_{\varepsilon,n_\varepsilon}}(\cdot))+3\varepsilon,$$
i.e. the last claim.
\hfill $\square$


\end{document}